\documentclass[reqno,oneside,11pt]{amsart}
\usepackage[foot]{amsaddr}
\usepackage[T1]{fontenc}
\usepackage[utf8]{inputenc}
\usepackage[top=1.5cm,bottom=2cm,right=2cm,left=2cm]{geometry}
\usepackage{amsmath,amsfonts,amssymb,amscd,amsthm}
\usepackage[usenames,dvipsnames,svgnames]{xcolor}
\usepackage{doi}
\usepackage{caption}
\usepackage{subcaption}
\usepackage{autonum} % Only enumerate what is referred

\usepackage{mathtools}
\usepackage[hang]{footmisc}
\usepackage{newtxtext}
\usepackage{newtxmath}
\usepackage{acronym}

\usepackage{booktabs}

\usepackage{tikz}
\definecolor{myellow}{RGB}{255,230,128}
\definecolor{gray20}{RGB}{204,204,204}
\definecolor{mygray}{RGB}{204,204,204}
\definecolor{mygreen}{RGB}{138,203,95}
\definecolor{myblue}{RGB}{77,151,214}

\usepackage[textsize=small]{todonotes}

\usepackage{soul}

\usepackage[ruled, lined, longend, linesnumbered]{algorithm2e}
\usepackage{xcolor}

\SetAlFnt{\footnotesize}
\SetArgSty{textnormal}

\SetAlCapSty{xAlCapSty}

\SetCommentSty{xCommentSty}

\SetNlSty{mynlfont}{}{} 

\LinesNumbered

\SetSideCommentRight

\DontPrintSemicolon

\RestyleAlgo{algoruled}

\SetInd{0.5em}{0.5em}

\SetAlgoCaptionLayout{SetSmall}

\newtheorem{theorem}{Theorem}[section]

\newtheorem{assumption}[theorem]{Assumption}

\tikzstyle{fig-ph}=[draw,minimum width=\textwidth, minimum height=\textwidth,text width=0.9\textwidth,color=red]

\tikzstyle{fig-ph-r}=[draw,minimum width=\textwidth, minimum height=1.4\textwidth,text width=0.9\textwidth,color=red]
\tikzstyle{fig-ph-rl}=[draw,minimum width=\textwidth, minimum height=0.5\textwidth,text width=0.9\textwidth,color=red]

\acrodef{dof}[DOF]{Degree Of Freedom}

\acrodef{UQ}[UQ]{Uncertainty quantification}
\acrodef{PDE}[PDE]{partial differential equation}
\acrodef{PDF}[PDF]{probability density function}
\acrodef{PRNG}[PRNG]{pseudo random number generator}
\acrodef{QoI}[QoI]{quantity of interest}
\acrodef{MC}[MC]{Monte Carlo}
\acrodef{MLMC}[MLMC]{multilevel Monte Carlo}
\acrodef{EMLMC}[EMLMC]{embedded \ac{MLMC}}
\acrodef{XFEM}[XFEM]{extended finite element method}
\acrodef{X-SFEM}[XS-FEM]{extended stochastic finite element method}
\acrodef{agfem}[AgFEM]{aggregated finite element method}
\acrodef{agfe}[AgFE]{aggregated finite element}

\acrodef{FE}[FE]{finite element}
\acrodef{CG}[CG]{conjugate gradient}

\DeclareTextFontCommand{\vbt}{\ttfamily\hyphenchar\font=45\relax}
\def\grad{{\boldsymbol{\nabla}}}
\def\diff{{k}}
\def\bm#1{\boldsymbol{#1}}
\def\x{{\boldsymbol{x}}}
\def\y{{\boldsymbol{y}}}
\def\X{{\boldsymbol{\xi{}}}}
\def\w{{\omega}}

\def\D{{\mathcal{D}}}
\def\M{{\mathcal{M}}}

\def\TN{{\mathcal{TN}}}

\def\F{{\mathcal{F}}}
 
\def\e{{\mathcal{E}}}

\def\U{{\mathcal{U}}}
\def\S{{\mathcal{S}}}
\def\P{{\mathcal{P}}}

\font\doble=msbm10 scaled\magstep1

\newcommand\R{\hbox{\doble R}}
\newcommand\E{\hbox{\doble E}}
\newcommand\V{\hbox{\doble V}}

\def\Fig#1{Fig.~(\ref{#1})}
\def\Eq#1{Eq.~(\ref{#1})}

\def\Ass#1{{\bf Assumption (\ref{#1})}}

\def\Sec#1{Section \ref{#1}}

\graphicspath{{figures/}}

\begin{document}

\title[EMLMC for UQ in random domains]{Embedded multilevel Monte Carlo for uncertainty quantification in random domains}

\author[S. Badia]{Santiago Badia$^{1,2}$}
\author[J. Hampton]{Jerrad Hampton$^{2}$}
\author[J. Principe]{Javier Principe$^{2,3,*}$}

\thanks{\null\\
$^1$ School of Mathematics, Monash University, Clayton, Victoria, 3800, Australia.\\
$^2$ Centre Internacional de M\`etodes Num\`erics en Enginyeria, Esteve Terrades 5, E-08860 Castelldefels, Spain.\\
$^3$ Universitat Polit\`ecnica de Catalunya, Campus Diagonal Bes\`os, Av. Eduard Maristany 16, Edifici A (EEBE), 08019, Barcelona, Spain\\
$^*$ Corresponding author.\\
E-mails: {\tt santiago.badia@monash.edu} (SB), {\tt jhampton@cimne.upc.edu} (JH), {\tt principe@cimne.upc.edu} (JP)}

\thanks{This research was supported by the European Union's Horizon 2020 research and innovation programme under the ExaQUte project, with grant agreement No 800898. JH has received funding from the European Union's Horizon 2020 research and innovation programme under the Marie Sk\l{}odowska-Curie grant agreement No 712949 (TECNIOspring PLUS) and from the Agency for Business Competitiveness of the Government of Catalonia.}
\date{\today}

\begin{abstract}
  The multilevel Monte Carlo (MLMC) method has proven to be an effective variance-reduction statistical method for Uncertainty quantification in PDE models. It combines approximations at different levels of accuracy using a hierarchy of meshes in a similar way as multigrid. The generation of body-fitted mesh hierarchies is only possible for simple geometries. On top of that, MLMC for random domains involves the generation of a mesh for every sample. Instead, here we consider the use of embedded methods which make use of simple background meshes of an artificial domain (a bounding-box) for which it is easy to define a mesh hierarchy, thus eliminating the need of body-fitted unstructured meshes, but can produce ill-conditioned discrete problems. To avoid this complication, we consider the recent aggregated finite element method (AgFEM). In particular, we design an embedded MLMC framework for (geometrically and topologically) random domains implicitly defined through a random level-set function, which makes use of a set of hierarchical background meshes and the AgFEM. Performance predictions from existing theory are verified statistically in three numerical experiments, namely the solution of the Poisson equation on a circular domain of random radius, the solution of the Poisson equation on a topologically identical but more complex domain, and the solution of a heat-transfer problem in a domain that has geometric and topological uncertainties. Finally, the use of AgFE is statistically demonstrated to be crucial for complex and uncertain geometries in terms of robustness and computational cost.
\end{abstract} %The abstract summarizes key findings in the paper and should be a paragraph no more than 250 words.

 %up to 10 key terms

\maketitle

%\noindent{\bf 2010 Mathematics Subject Classification:} 65N12, 65N15, 65N30


%\tableofcontents

\section{Introduction}
\label{sec:intro}
\ac{UQ} requires the solution of a stochastic \ac{PDE} with random data. Several methods for solving stochastic \acp{PDE}, such as stochastic Galerkin \cite{lemaitre2010spectral,ghanem1991stochastic} or stochastic collocation \cite{Babuska2010}, are based on a standard approximation in space like \acp{FE} or finite volumes, and different types of polynomial expansions in the stochastic space \cite{Xiu2003}. Although very powerful for some particular problems, most of these techniques suffer the so-called `curse of dimensionality', i.e. a dramatic increase of computational cost with the number of stochastic variables. Besides, some of these methods are intrusive because given a code that can be used to solve a deterministic problem, it needs to be modified to solve a stochastic one.

In contrast, the classical \ac{MC} method does not present these drawbacks. On the one hand it is a sampling method, i.e., it only requires the repeated evaluation of a deterministic model, and it can be implemented in a non-intrusive way. On the other hand it is known to convergence to the exact statistics of the solution as the number of input samples tends to infinity, independently of the dimensionality of the stochastic space and mostly independently of the physics of the problem under consideration, as long as some moments of the \ac{QoI} are bounded. However, the number of samples required to achieve statistical convergence combined with the complexity of the computational model required to have enough spatial and/or temporal accuracy can make it prohibitively expensive.

To address this difficulty, the \ac{MLMC} method was introduced in \cite{Kebaier2005} (with two levels) and in \cite{Giles2008} (with several levels) and applied to solve stochastic \acp{PDE} with random input parameters in \cite{Giles2012,Mishra2012,Mishra2012a,Pisaroni2017}. The \ac{MLMC} method exploits  a hierarchy of  discretizations of the underlying differential problem (using a hierarchy of meshes). The method computes expectations on the difference between solutions at two consecutive lvels of accuracy, reducing the variance of the problem at hand and thus decreasing the number of samples required at finer discretization levels. By doing this, most of the computational cost is transferred to coarse discretization levels, where a large number of samples is used to control the statistical error of the \ac{MLMC} estimator, whereas only few samples are used on the finest (and most costly) levels, to control the discretization error.

Depending on the application, different types of uncertain data can be considered, like material properties (\ac{PDE} parameters), boundary and/or initial conditions, and the geometry/topology of the domain. The \ac{MLMC} method can be used out-of-the-box in most cases but it requires a hierarchy of discretizations as input, which can be difficult to generate when dealing with complex geometries and unstructured body-fitted meshes. For this reason \ac{MLMC} methods are barely used to perform \ac{UQ} when complex geometries are considered. This is the first limitation we aim to address herein, utilizing embedded methods with \ac{MLMC} in a method referred to as \ac{EMLMC}.

On top of that, leaving simple academic constructions aside, geometry is hardly known exactly and surfaces are rarely smooth when looking at sufficiently small scale. Variation between different samples of natural materials, e.g. porous rocks, and cost limits of diminishing manufacturing tolerances of fabrication processes make geometry an epistemic uncertainty in general. Geometric uncertainties can have a strong impact on the solution of a \ac{PDE}, e.g. they trigger boundary layer separation in fluids, and these effects are important in many fields, e.g. geology \cite{Deffeyes1982}, tribology \cite{hutchings2017tribology} and lubrication \cite{hamrock2004fundamentals}, nano-technology \cite{Nosonovsky2008} and biology \cite{Park2012}. The first attempts to predict the impact of roughness was through deterministic parameterizations, like periodic indentations \cite{Taylor1971,Richardson1971}, sinusoidal corrugation \cite{Fyrillas2001} and fractal representations, e.g. using the well-known Von Koch snowflake curve \cite{Cajueiro1999,Blyth2003}. After that, stochastic representations of surface roughness have been adopted frequently, which permit to incorporate the lack of information or measurement errors. Examples of these descriptions include random fields to represent the thickness in shell structures (see \cite{Stefanou2004} and the references therein), random fractal representations \cite{Sarkar1993} and random fields with a given spatial correlation structure to represent derivations from a nominal geometry \cite{Xiu2006,Tartakovsky2006,Mohan2011874,Zayernouri2013}. In engineering design, geometry is exactly defined using computational tools and in this case the uncertainty is rooted at manufacturing process, as mentioned. Uncertainty also arises when measuring, e.g  when the geometry is determined from digital images, sourced at image scanning and resolution. Using pixelated images for computing may result in significant errors, even in the limit of small resolution \cite{Babuska2001,Babuska2003}. However, it is possible to recover a definition of the geometry using level set methods \cite{sethian1999level} and, in fact, this can be done statistically, e.g. using polynomial chaos expansions, to define random level set representations \cite{Stefanou2009}. Curiously, level set and phase field representations can also be obtained from  noisy images solving stochastic \acp{PDE} \cite{Patz2013}.

The are three approaches to deal with geometric uncertainty. The first one is based on domain mappings, the second on perturbation methods and the third one on the use of fictitious domain methods, also known as embedded or immersed methods. The first approach started with the stochastic mapping concept that was introduced in \cite{Xiu2006} to transform a (deterministic or stochastic) problem in a random domain to a stochastic problem in a deterministic domain. This approach was applied to the stochastic analysis of roughness in flow problems in \cite{Tartakovsky2006}. In \cite{Mohan2011874} the main idea is to fix the mesh connectivity and change the coordinates of the nodes using deformation strategies that incorporate uncertainty through the solution of auxiliary \ac{PDE} with random boundary conditions. In \cite{Castrillon-Candas2016} a domain mapping approach is used to prove the convergence of a stochastic collocation method based on Smolyak grids. Because the construction of globally smooth mappings is difficult, a piecewise smooth domain mapping approach is proposed in \cite{Chaudhry20181127}. The second approach was introduced in \cite{Sarkar1993} where an analytical perturbation method was used to find solutions of the Laplace equation. A perturbation method developed using shape calculus (a Taylor expansion using shape derivatives and its linearization around a nominal domain) was proposed in \cite{Harbrecht2008,Harbrecht2013,Harbrecht2016}. Because shape derivatives are difficult to compute, a perturbation approach based on approximated stochastic boundary conditions (of Robin type) in \cite{Dambrine2017943,Dambrine2016921}. Although it is not based on a formal perturbation expansion, the stochastic smoothed profile method of \cite{Zayernouri2013} can be included in this category. This method transforms a random domain into a random force-term in a deterministic domain through the introduction of a smoothly spreading interface layer to represent rough boundaries and its error grows with the layer thickness \cite{Luo2008}. By definition, these approaches are not applicable to complex geometrically and even topologically random domains.

The use of embedded methods, which rely on a fixed background grid and a discrete approximation of the geometry using a level set function, for the solution of \acp{PDE} in stochastic domains combined with polynomial chaos was proposed in \cite{Canuto2007} and applied to fluid dynamics problems in \cite{Parussini2010}. Among the several possible immersed methods, those in \cite{Canuto2007,Parussini2010} impose boundary conditions through Lagrange multipliers defined in a \ac{FE} space built from a boundary discretization (independent of the domain discretization). This construction results in a saddle point linear system for which specific preconditioners have been recently developed \cite{Gordon2014}. Alternatively, the use of the so-called \ac{XFEM} method \cite{Moes1999,sukumar_modeling_2001} for dealing with embedded stochastic geometries using polynomial chaos, was proposed in \cite{Nouy2007,Nouy2008,Nouy20101312,Nouy20113066} and named \ac{X-SFEM}. It was lately applied to heat transfer \cite{Lang20131031} and structural problems \cite{Schoefs2016}. The lack of robustness of embedded \ac{FE} methods is not addressed in these works. 

The use of \ac{MLMC} in complex random domains using body-fitted meshes would require to generate an unstructured mesh (and its corresponding mesh hierarchy) per sample, which is absolutely impractical. For this reason, and the previously described problems related to the generation of mesh hierarchies in general, we also consider embedded \ac{FE} methods in this work. However, the naive use of methods exhibit some problems. The most salient one is the so-called small cut cell problem. The intersection of cells with the domain surface cannot be controlled and one can end up with cells for which the ratio between the cell volume and its portion in the domain interior is arbitrarily large. It is well-known that the condition number of the resulting linear system blows up with these ratios. As a result, plain embedded \ac{FE} methods are not robust. Different remedies have been proposed so far in the literature. For conforming \ac{FE} methods in which trace continuity must be enforced between cells, one can consider methods that introduce artificial dissipation (stabilization) in order to weakly enforce zero jump of derivatives (up to the \ac{FE} order being used) of jumps on facets that belong to cut cells (see \cite{burman_cutfem:_2015} for more details). Another approach is to create aggregated meshes such that every aggregate contains at least one interior cell. Whereas the generation of discontinuous \ac{FE} spaces on aggregated meshes is straightforward (enforcing the local \ac{FE} space to be the same polynomial space as in a standard cell), it is far more complicated for conforming \acp{FE}. Using such a naive approach for conforming methods would generate global constraints and too much \emph{rigidization}, which would affect the convergence properties of the method, its implementation, and the sparsity of the resulting linear system; we note that this is in fact the \emph{strong} version of the weakly enforced constraints in \cite{burman_cutfem:_2015}.  The aggregated \ac{FE} method has been proposed in \cite{Badia2017a} to enable aggregation techniques for conforming methods. It has three key features: (a) the constraints are cell-local and their implementation in \ac{FE} codes is easy; (b) it keeps the convergence order of standard body-fitted \ac{FE} spaces; (c) it does not perturb the Galerkin formulation with any kind of artificial dissipation and does not involve the computation of (possibly) high-order derivatives on the element boundaries.

In this work, we propose the \ac{EMLMC} framework. It combines the \ac{MLMC} method on a hierarchy of background meshes (which can be simple structured meshes) and the \ac{agfem} at every mesh level to capture every sample of the random domain. In particular, we consider random domains implicitly defined through random level-set functions and the marching tetrahedra scheme to approximate the random domain at every mesh level. The resulting method extends the applicability of \ac{MLMC} to complex random domains and it is robust. We state the model problem in Section \ref{sec:problem} and we describe the \ac{MLMC} method in Section \ref{sec:mc_and_mlmc}. The construction of the \ac{agfem} discretization for each sample is then described in Section \ref{sec:agg}. Finally, a set of numerical experiments are described in Section \ref{sec:examples} and we draw some concluding remarks in Section \ref{sec:conclusions}.

\def\artdom{\mathcal{B}}
\section{Problem formulation}
\label{sec:problem}
As a model problem we consider the following elliptic stochastic problem. Given an oriented manifold $\M(\w)\subset \mathbb{R}^d$ and its corresponding interior domain $\D(\w)$, find $u(\x,\w)$ such that
\begin{align}
  -\grad \cdot \left( \diff \grad u \right) = f  \ \hbox{ in } \ \D(\w), \qquad
  u = u_0 \ \hbox{ on } \, \M(\w), \label{eq:model}
\end{align}
almost surely (a.s.), i.e., for almost all $\w \in \Omega $, which denotes the uncertainty, described by a complete probability space $(\Omega,\F,P)$. Although stochastic coefficients, i.e. the diffusion $\diff = \diff(\x,\w)$, forcing $f=f(\x,\w)$, and boundary condition $u_0 = u_0(\x,\w)$ could be considered random fields too, we assume they are deterministic, as it is usual in the literature when stochastic domains are considered \cite{Chaudhry20181127,Mohan2011874,Xiu2006,Harbrecht2008,Harbrecht2013,Harbrecht2016,Dambrine2017943,Dambrine2016921}. Let us assume that all realizations $\M(\w)$ (and $\D(\w)$) are bounded. We can define a bounded artificial domain $\artdom$ that contains all possible realizations of $\M(\w)$, i.e. $\D(\w) \subset \artdom, \, \forall \w \in \Omega$. We also assume that $\diff$ and $f$ are defined in $\artdom$, independently of $\w \in \Omega$. With the random solution of this problem at hand we aim to compute $\E\left(Q(u)\right)$ where $Q$ is a deterministic \ac{QoI}, e.g. an integral of $u$ on a subregion or surface, and $\E$ the expectation.

Although it is not strictly necessary (at least in deterministic domains \cite{Charrier2013,Teckentrup2013a}), uniform ellipticity is often assumed to guarantee the well-posedness of \Eq{eq:model} (see \cite{Barth2011,Babuska2010,Chaudhry20181127}). Therefore, we assume $\diff$ uniformly bounded from below, i.e. there exists $\diff_0$ such that $\diff(\x) \ge \diff_0, \, \forall \x \in \artdom$ (and $\forall \w \in \Omega$ if a stochastic diffusion $\diff(\x,\w)$ is considered). Under these assumptions, the bilinear form associated to the weak form of \eqref{eq:model} is bounded and coercive and the Lax-Milgram lemma guarantees a solution for any $\w \in \Omega$ uniformly bounded by $ \|f\|_{L^2( \artdom)} $ \cite{Chaudhry20181127}.

\def\S{\mathcal{S}}
\def\M{\mathcal{M}}
\def\T{\mathcal{T}}
\def\lev{{\ell}}

\def\Tcut{{\T}^\M}
\def\Tact{{\T}^\D}
\def\Dact{\widetilde{\D}}
\def\X{\mathcal{V}}
\def\Xact{\widetilde{\X}}

\def\I{\mathcal{I}}

\section{Monte Carlo and its multilevel extension}
\label{sec:mc_and_mlmc}
\subsection{The Monte Carlo method}
\label{sec:mc}
The aim is to compute the expected value $\E(Q)$ of the \ac{QoI} $Q(u)$. Using a discrete approximation in a mesh $\T_h$ of size $h$ having $M \propto h^{-d} $ degrees of freedom, we actually compute a discrete approximation to the \ac{QoI} using $u_h$, the discrete \ac{FE} solution of the continuous problem at hand. Let us denote $Q_h=Q(u_h)$. The standard MC algorithm computes the approximation of its expected value as the average
\begin{equation}
\E(Q) \approx \overline{Q}_h = \frac{1}{N} \sum_{i=1}^N Q_h^i, \label{eq:mc_average}
\end{equation}
where $Q_h^i$ are obtained evaluating the \ac{QoI} using the $N$ realizations (samples) $u_h^i$ of the stochastic solution $u_h(\x,\w)$ on the given mesh $\T_h$. The mean square error of this approximation is
\begin{equation}
e^2(\overline{Q}_h) = \E\left[(\overline{Q}_h - \E(Q))^2\right] =  (\E(Q_h -Q))^2 + \V(\overline{Q}_h)
\label{eq:mc_error}
\end{equation}
where $\V$ is the variance. As it is well known, the first term is the (squared) discretization or bias error and its reduction requires refining the mesh, whereas the second one is the statistical error that decays as
\begin{equation}
\V(\overline{Q}_h) =  N^{-1} \V(Q_h),
\label{eq:mc_variance}
\end{equation}
so its reduction requires increasing the number of samples. Denoting the complexity of evaluating $Q^i_h$ as $C_h$, it is possible to estimate the total computational cost of the MC algorithm, $C_{\rm MC}$, under the following assumptions.
\begin{assumption}
  \label{as:1}
  There exist $\alpha$ and $c_{\alpha}$ such that
  $$|\E(Q_h -Q)| \leq c_{\alpha}h^{\alpha} .$$
\end{assumption}
\begin{assumption}
  \label{as:2}
  There exist $\gamma$ and $c_{\gamma}$ such that
  $$C_h \leq c_{\gamma} h^{-\gamma} .$$
\end{assumption}

By \Ass{as:1} the mesh size required to achieve a discretization error smaller than $\varepsilon$ is $h \preceq  \varepsilon^{1/\alpha}$ (hereafter $a \preceq b$ means that there is a constant $c$, independent of $h$, such that $a \le c b$). 
On the other hand, by \Eq{eq:mc_variance}, the number of samples required to achieve a statistical error $ \left[\V(\overline{Q}_h)\right]^{1/2} \leq  \varepsilon $ is $N \propto \varepsilon^{-2}$, where here $a\propto b$ means that there exist constants $c$ and $C$ such that $c\varepsilon^{-2} \le N \le C\varepsilon^{-2}$. Then
\begin{equation}
C_{\rm MC} = N C_h \preceq  N h^{-\gamma}  \preceq \varepsilon^{-(2+\gamma/\alpha)}
\label{eq:mc_complex}
\end{equation}
Note that the factor $\varepsilon^{-\gamma/\alpha} $ comes from the discretization error and represents the (asymptotic) cost of one sample.
For the model problem of \Sec{sec:problem}, \Ass{as:1} holds with $\alpha=2$ \cite{Teckentrup2013a}. Now the constant $\gamma$ depends on $d$ and the type of solver. For naive Gaussian elimination $\gamma=3d$, whereas $\gamma=d$ for the optimal multigrid method. In between, for sparse linear solvers $\gamma=3(d-1)$ (for $d=2,3$) whereas for the (unpreconditioned) \ac{CG} algorithm $\gamma=d+1$. In the last case, for example, the computational cost scales as  $\varepsilon^{-4}$ when $d=3$. 

\subsection{The Multilevel Monte Carlo method}
\label{sec:mlmc}

The ac{MLMC} algorithm exploits the linearity of the expectation using a hierarchy of $L+1$ meshes $\T_0,\T_1,...,\T_L$ of sizes $h_0>h_1>...>h_L$. Although there are other options \cite{Haji-Ali2016}, we assume $h_l=h_0s^{-l}$ (i.e. each mesh in the hierarchy is obtained by uniformly dividing each cell into $s^d$ subcells). Performing $\{N_\lev\}_{\lev=0,...,L}$ simulations for different values of the random parameters $\w$ on $\{\T_\lev\}_{\lev=0,...,L}$, the expectation on the coarse grid is corrected using the whole hierarchy as
\begin{equation}
\E(Q_L) = \E(Q_0) + \sum_{\lev=1}^L \E(Q_\lev-Q_{\lev-1}) \approx \sum_{\lev=0}^L \overline{Y_\lev}:=\widetilde{Q}_L \label{eq:mlmc_average}
\end{equation}
where $Q_\lev=Q(u_{\lev})$ is the approximation of the \ac{QoI} computed using the \ac{FE} solution $u_\lev$ for the mesh $\T_\lev$ and $Y_\lev=Q_\lev-Q_{\lev-1}$ for $\lev=1,...,L$ and $Y_0=Q_0$.
The total error of this approximation is now given by
\begin{equation}
  e^2( \widetilde{Q}_L ) = \E\left[( \widetilde{Q}_L - \E(Q) )^2\right] =  (\E(Q_L -Q))^2 + \V(\widetilde{Q}_L)
\label{eq:mlmc_error}
\end{equation}
The first term, the (squared) discretization error, is the same as in the \ac{MC} method and it requires mesh $\T_L$ to be fine enough. For this term to be smaller than, e.g. $\varepsilon^2/2$, \Ass{as:1} requires $c_{\alpha} h_L \leq 2^{-1/2} \varepsilon^{1/\alpha}$, or, equivalently,
\begin{equation}
  L = \lceil {\frac{1}{\alpha} \log_s \left(\sqrt{2}c_{\alpha} \varepsilon^{-1}\right)} \rceil
 \label{eq:mlmc_num_levels}
\end{equation}
where, as usual, $\lceil x \rceil$ denote the unique integer $n$ satisfying $x<n<x+1$.  The second term \Eq{eq:mlmc_error} is the statistical error, given by $\V(\widetilde{Q}_L) = \sum_{\lev=0}^{L} N_\lev^{-1} \V(Y_\lev)$, is very different from the one in the \ac{MC} method thanks to the following assumption.
\begin{assumption}
  \label{as:3}
  For any $\lev=0,....L$, there exist $\beta$ and $c_{\beta}$ such that
  $$\V(Y_\lev) \leq c_{\beta}h_\lev^{\beta}.$$
\end{assumption}
Thanks to \Ass{as:3}, $\V(Y_\lev)$ tends to zero with $h_\lev$ and in this sense, the \ac{MLMC} method can be understood as a variance reduction method. Under \Ass{as:1}, \Ass{as:2} and \Ass{as:3}, it is possible to explicitly bound the total computational cost and to minimize it to obtain the optimal number of samples to be taken on each level \cite{Giles2008}.

The number of samples on each level required to have a (squared) statistical error smaller than $\varepsilon^2/2$ is given by
\begin{equation}
  N_\lev =   \lceil {2 \varepsilon^{-2} \sqrt{\frac{\V(Y_\lev)}{C_\lev}}  \sum_{i=0}^{L} \sqrt{\V(Y_i)C_i}} \rceil,
\label{eq:mlmc_num_samples}
\end{equation}
where $C_\lev$ is the complexity (computational cost) of evaluating $Q_\lev$, i.e. taking one sample at level $\lev$. \Eq{eq:mlmc_num_samples} can also be written as
\begin{equation}
  N_{\lev} = \lceil {s^{\Gamma (L-\lev)} N_{L}} \rceil,
  \label{eq:num_samples_ratio}
\end{equation}
where
$$\Gamma=\frac{1}{2}(\gamma+\beta).$$

The computational complexity of the \ac{MLMC} algorithm is given by $C_{\rm MLMC} = \sum_{\lev=0}^{L} N_\lev C_\lev$,  and it can be bounded using \Ass{as:2} and \Ass{as:3} as
\begin{equation}
 C_{\rm MLMC} \preceq \varepsilon^{-2-\max(0,\frac{\gamma-\beta}{\alpha})},
\label{eq:mlmc_complexity_bound}
\end{equation}
with an additional logarithmic factor when $\gamma=\beta$.

There are three different situations depending on whether the (computational cost required to reduce the) statistical error dominates, is of the same order of, or is dominated by the discretization error. In the first case, $\beta>\gamma$ and because the statistical error is eliminated using coarse grids, the complexity of the \ac{MLMC} algorithm is independent of the fine mesh size. In the last one $\beta<\gamma$ and the complexity includes a factor $\varepsilon^{-\gamma/\alpha} $ that comes from the discretization error, which also appears in the \ac{MC} complexity estimate \eqref{eq:mc_complex}. In applications of \acp{PDE} with random coefficients, it is typical to have $\beta=2\alpha$ and in this case, $C_{\rm MLMC} \leq \varepsilon^{-\gamma/\alpha} $ which is the complexity of one sample at the finest grid, that is, the \ac{MLMC} algorithm does not  \emph{see} the sampling cost \cite{Pisaroni2017}. Thus, the relation between parameters $\beta$ and $\gamma$ determines whether the effort should be concentrated on coarse or fine meshes.

From the previous discussion the crucial role of the constants  $\alpha$, $\beta$, $\gamma$, $c_{\alpha}$, $c_{\beta}$ and $c_{\gamma}$ is apparent. In general, these constants are unknown (they are problem and solver dependent) but their estimation is required to define the optimal number of levels and samples in \ac{MLMC}. The so-called adaptive \ac{MLMC} algorithms have been developed to face this challenge by estimating these constants with extrapolations and/or statistical inference \cite{Giles2008,Cliffe2011,Collier2015}. For example, the variances in \Eq{eq:mlmc_num_samples} can be estimated (a posteriori) using the sample variance
\begin{equation}
  \V(Y_\lev) \approx V(Y_\lev) = \frac{1}{N_\lev} \sum_{i=1}^{N_\lev} (Y^i_\lev -  \overline{Y_\lev})^2.
  \label{eq:sample_variance}
\end{equation}
We do not follow this strategy herein, but we exploit the knowledge of the problem at hand.
For the problem of \Sec{sec:problem} $\beta=4$ \cite{Teckentrup2013a} and using a \ac{CG} algorithm $\gamma=d+1$, so the statistical error dominates the computational cost and the optimal relation between samples at the hierarchy are given by \Eq{eq:num_samples_ratio} with $\Gamma=4$ for $d=3$ and $\Gamma=3.5$ for $d=2$.

Observe that  \Ass{as:1} and \Ass{as:3} (redefining $c_{\alpha}$ to include $h_0^{-\alpha}$ and $c_{\beta}$ to include $h_0^{-\beta}$ ) can be written as
\begin{align}
  |\E(Q_\lev-Q)| \leq  c_{\alpha} s^{-\alpha \lev}, \qquad \V(Y_\lev) \leq  c_{\beta} s^{-\beta \lev}.\label{eq:error_decay} 
\end{align} 
The corresponding sample approximations $|Q-\widetilde{Q}_\lev|$ (see \Eq{eq:mlmc_average}) and $V(Y_\lev)$ (see \Eq{eq:sample_variance}) statistically converge to these values and thus, with a large enough number of samples, one can check spatial convergence (i.e., with respect to $\ell$) for these quantities too.
In the numerical examples of \Sec{sec:examples} $s=2$ and therefore, the (logarithmic) convergence rate with respect to the level is $2\log(2)$ and $4\log(2)$ respectively.

\section{The aggregated FE method}
\label{sec:agg}
The mesh hierarchy required by the \ac{MLMC} algorithm is built from $\T_0$, an initial shape-regular partition of $\artdom$, generating $\T_{\lev+1}$ by uniform refinement of $\T_\lev$. 
In our implementation, we consider (a bounding box) $\artdom := \prod_{i=1}^d[x^i_m,x^i_M]$ and $\T_0$ a uniform Cartesian mesh. At every level $\lev$, we consider a piecewise linear approximation $\M_\lev(\w)$ of the manifold $\M(\w)$. When $\M(\w)$ is already a polytopal surface mesh, e.g. an STL mesh, no approximation is required at this step. However, if the resolution of $\M(\w)$ is very fine compared to the coarser mesh size, it would involve a very costly numerical integration at coarser levels. A case of practical interest that will be addressed in this work is when $\M(\w)$ is represented implicitly with a level-set function $\phi(\x)$. In this case, we construct a $\mathcal{C}^0$ Lagrangian \ac{FE} space of arbitrary order on the grid $\T_\lev$, which is represented with ${\Xact}(\T_\lev)$. In this case, one can define $\M_\lev(\w)$ as the result of applying the marching tetrahedra algorithm \cite{Doi91} to the interpolation $\phi_\lev(\x)$ of the level-set function $\phi(\x)$ onto $\Xact(\T_\lev)$. When considering quad meshes, one can simply use the decomposition of an n-cube into n-tetrahedra before applying the marching tetrahedra algorithm.
 
\newcommand{\closure}[2][3] {
{}\mkern#1mu\overline{\mkern-#1mu#2}}
Given the oriented manifold $\M_\lev(\w)$ (and its corresponding interior, the domain $\D_\lev(\w)$) and the background mesh $\T_\lev$, due to the piecewise linear nature of the manifold, it is possible to compute exactly the intersection of each cell $K \in \T_\lev$ with $\M_\lev(\w)$ and $\D_\lev(\w)$. We represent with $\Tact_\lev(\w)$ (resp. $\Tcut_\lev(\w)$) the set of active (resp., cut) cells in $\T_\lev$ that intersect $\D_\lev(\w)$ (resp., $\M_\lev(\w)$); $\Tact_\lev(\w) \setminus \Tcut_\lev(\w)$ is the set of interior cells in $\D_\lev(\w)$. 
Numerical integration over cells $K \in \Tcut_\ell(\w)$  can be carried out by generating a simplicial mesh $\I^\D_K$ of $K \cap \D_\lev(\w)$, e.g. using a Delaunay triangulation. 
Interior cell volume integration is straightforward, and we can simply take $\I_K^{\D} = K$. We proceed analogously to create surface meshes  for ${K} \cap \M_\lev(\w)$. We define the volume integration mesh as the union of the interior cells and the integration meshes for cut cells $ \I_\lev^\D(\w) := \bigcup_{K \in \Tact_\ell(\w)} I^\D_K$; analogously for the surface mesh $\I_\lev^\M(\w)$. These integration meshes always exist and are cell-wise local, i.e., no conformity is required among cells.
We note that $\I_\lev^\D(\w)$ is a body-fitted mesh of $\D_\ell(\w)$ and $\I_\lev^\M(\w)$ is a surface mesh of $\M_\lev(\w)$. However, the shape regularity of the meshes is not relevant since they are only used for integration purposes. For higher (than one) order approximations of $\M(\w)$, we can consider the same approach described above supplemented with an additional degree elevation of the linear cell-wise integration meshes. These techniques impose some constraints on the ratio between surface curvature and mesh resolution and are not being used in Section \ref{sec:examples}, since we restrict ourselves to first order \acp{FE}.

\def\gradient{\boldsymbol{\nabla}}
In this work, $H^1$-conforming Lagrangian space of arbitrary order are considered for the discretization of the \ac{PDE} problem at hand. We denote with $\Xact_\ell(\Tact,\w)$ the \ac{FE} space on the grid $\Tact_\lev(\w)$. The space $\Tact_\lev(\w)$ is defined on $\Dact(\w) \supseteq \D(\w)$, where $\Dact(\w)$ is the union of all cells in $\Tact(\w)$. Thus, it is not a body-fitted mesh of $\D(\w)$ and embedded \ac{FE} techniques must be used to properly capture the geometry. Let us consider the weak form of \eqref{eq:model} in which the essential boundary conditions are imposed in a weak form using Nitsche's method. First, we define the cell-wise forms
\def\A{\mathcal{A}}
\def\n{\boldsymbol{n}}
\def\cell{K}
\begin{align}
  &  \A^\lev_\cell(u,v) := \int_{\cell \cap \D_\lev(\w)} \gradient u \cdot \gradient v \mathrm{d}V + \int_{\cell \cap \M_\lev(\w)} (\tau_\cell u v - v ( \n \cdot \gradient u) - u (\n \cdot \gradient v) ) \mathrm{d}S, \\
  & {b^\lev_\cell(v)} := \int_{\cell \cap \D_\lev(\w)} f v \mathrm{d}V + \int_{\cell \cap \M_\lev(\w)} (\tau_K u_0 - (\n \cdot \gradient v) u_0  ),
\end{align}
where $\tau_\cell \approx \mathcal{O}(h_\lev^{-2})$ is a positive parameter that must be large enough to ensure coercivity of the bilinear form.  The global form $\A^\lev: H_0^1(\D_\lev(\w)) \rightarrow H^{-1}(\D_\lev(\w))$  and right-hand side term $b^\lev \in H^{-1}(\D_\ell(\w))$ are stated as the sum of the element contributions, i.e.,
\begin{align}\label{eq:bil-form}
  \A^\lev(u,v) := \sum_{\cell \in \Tact_\lev(\w)} \A^\lev_\cell(u,v),  \quad
  b^\lev{ (v)} := \sum_{\cell \in \Tact_\lev(\w)} b^\lev_\cell(v).
\end{align}
The use of embedded \ac{FE} methods has a dramatic impact on the condition number of the resulting linear systems, known as the \emph{small cut cell} problem. Given that the mesh $\T_\lev$ is shape regular, one can define its characteristic cell size $h_\lev$. It is well-known in body-fitted \ac{FE} methods that the condition number of the linear system that results from the discretization of second order elliptic \acp{PDE} is  $\mathcal{O}(h_\ell^{-2})$. However, for unfitted \ac{FE} methods, it also deteriorates with the maximum of ${|K \cap \D_\lev(\w)|}^{-1}{|K|}$ among all cells in $K \in \Tact_\lev$. The portion of the cell in the physical domain cannot be controlled and it can be arbitrarily close to zero. As a result, embedded \ac{FE} methods can produce almost singular linear systems to be solved, and thus are not reliable. 

\def\Q{\mathcal{Q}}
\def\G{\mathcal{G}}
To fix the small cut cell problem issue, we consider the \ac{agfem} recently proposed in \cite{Badia2017a}. The idea is to start with the mesh $\Tact_\lev$ and generate an aggregated mesh $\G_\lev(\w)$ in which the cut cells with small volume ratio within the physical domain are aggregated to one interior cell (see \cite[Alg. 2.1]{Badia2017a}). Such aggregation is next used to define constraints over the standard \ac{FE} space $\Xact(\Tact_\lev,\w)$. We represent the resulting \ac{agfe} space with $\X(\G_\lev,\w) \subseteq \Xact(\Tact_\lev,\w)$. We refer the interested reader to \cite{Badia2017a} for a detailed exposition of the mesh aggregation algorithm, the computation of constraint, the implementation issues, and its numerical analysis, to \cite{Verdugo2019} for its parallel implementation, and to \cite{Badia2018b} for a mixed \ac{agfe} space for the Stokes problem. \ac{agfe} spaces are not affected by the small cut cell problem and the condition number of the resulting matrix is $\mathcal{O}(h_\lev^{-2})$ \cite{Badia2017a}. The numerical experiments in \cite{Verdugo2019} show that a standard parallel algebraic multigrid solver is effective to solve \ac{agfe} linear systems using default settings. The \ac{agfe} space approximation of \eqref{eq:model} at a given level $\lev$ reads as follows: find $u_\lev \in \X(\G_\lev,\w)$ such that 
\begin{align}
  \A^\lev(u_\lev,v_\lev) = b^\lev(v_\lev) \qquad \hbox{ for any } \, v_\lev \in \X(\G_\lev,\w).\label{eq:agfeprob}
\end{align}
We note that all these terms can be computed using the cell-wise integration triangulations commented above. The condition number of the matrix resulting from the \ac{agfem} discretization has been proven to be $\mathcal{O}(h_\lev^{-2})$ as in body-fitted methods (see \cite{Badia2017a}). All the steps required to apply \ac{agfem} in the frame of \ac{MLMC} are listed in Alg. \ref{alg:alg}.

\def\sample{\w^i}
\begin{algorithm}
  \caption{\ac{EMLMC}\label{alg:alg}}
\For{$\lev \in \{1,\ldots,L\}$}{
  Compute $\T_\lev$ as uniform refinement of $\T_{\lev-1}$\\
  Generate the \ac{FE} space $\Xact(\T_\lev)$
}
$\widetilde{Q}_L = 0 $ \\
\For{$\lev \in \{0,\ldots,L\}$}{
  $\overline{Y_l} = 0$\\
  \For{$i \in \{1,\ldots,N_\lev\}$}{
    Sample the random variables and get $\sample$ \\
    Set $k=\lev$ \\
    Interpolate the level-set function $\phi(\x,\sample)$ in the \ac{FE} space $\Xact(\T_k,\sample)$ to obtain $\phi_k(\x,\sample)$ \nllabel{lin:interp} \\
    Run the marching tetrahedra algorithm for $\phi_k(\x, \sample)$ and compute $\M_k(\sample)$ and $\D_k(\sample)$ \\
    Compute the intersection of cells in $\T_k(\sample)$ with $\M_k(\sample)$ and $\D_k(\sample)$ to obtain $\Tact_k(\sample)$ and $\Tcut_k(\sample)$ and the cell-local integration meshes $\I_k^{\M}(\sample)$ and $\I_k^{\D}(\sample)$ \\
    Generate the \ac{FE} space $\Xact(\Tact_k,\w^*)$ \\
    Run the mesh aggregation algorithm in \cite[Alg. 2.1]{Badia2017a} to compute the aggregated mesh $\G_k(\sample)$ \\
    Generate the \ac{agfe} space $\X(\G_k,\w)$ as in \cite[Sect. 3]{Badia2017a}\\
    Compute the solution $u_k(\w^*)$ of \eqref{eq:agfeprob}  \\
    Compute the $Q_k=Q(u_k)$  \nllabel{lin:qoi} \\
    $\overline{Y_l} = \overline{Y_l} + Q^i_k$\\
    \If{$\lev > 0$}{
      Set $k=\lev-1$ \\
      Repeat lines \ref{lin:interp} to \ref{lin:qoi}.\\
      $\overline{Y_l} = \overline{Y_l} - Q^i_k$\\
    }
  }
  $\widetilde{Q}_L = \widetilde{Q}_L + \overline{Y_l} / N_l$ \\
}
\end{algorithm}

\section{Numerical examples}
\label{sec:examples}

In this section we present three numerical examples aimed at illustrating the efficiency and robustness of the proposed approach. The first of them is a convergence test whose goal is to asses the actual convergence rates in the implementation. The second one is the solution of  \eqref{eq:model} in a complex random domain illustrating in a quantitative way the robustness of the proposed approach. Finally we solve \eqref{eq:model} in a situation where the randomness may induce changes in the topology. These numerical experiments have been carried out using \texttt{FEMPAR} v1.0.0~\cite{fempar,Badia2019}, an open source object-oriented scientific computing library for the simulation of complex multiphysics problems governed by \acp{PDE} at large scales that incorporates embedded \ac{FE} machinery \cite{Verdugo2019}. 

The implementation of \ac{EMLMC} was made in \texttt{FEMPAR} by developing a new software layer which permits the concurrent execution of a model, understood as the actual implementation of an algorithm that permits the solution of a \ac{PDE}. In our case, the model was developed using the numerical tools provided by  \texttt{FEMPAR} to perform the automatic generation of background meshes $\T_\lev$ of arbitrary size, the level-set description of the geometry and its interpolation, the construction of \ac{agfem} discretization, the numerical integration and the solution of the linear systems using domain decomposition preconditioners. The implementation of sampling methods made in this new module \cite{workflow} exploits three levels of parallelism, across levels, across samples and in the evaluation of each sample, e.g. employing a domain decomposition method, although the latter feature has not been exploited in the examples below (each sample is computed using one processor).

\subsection{A convergence test}
\label{subsec:Canuto}
In this section we consider an example with analytic solution proposed in~\cite{Canuto2007}.
The spatial domain $ \D(R)\subset \R^2$ is a circle centered at $(0.5,0.5)$ whose radius $R$ is a truncated Gaussian random variable with mean $\mu=0.3$, standard deviation $\sigma=0.025$, and restricted to $[a,b]$ with $a=0.2$ and $b=0.4$, which we denote by $\TN(a,b,\mu,\sigma)$. In this random domain we solve problem \Eq{eq:model} with $f=4$, $k=1$, and Dirichlet boundary conditions $u_0=u_{\rm ex}$ where 
\begin{equation}
  u_{\rm ex}(\x):=R^2-|\x_1 - 0.5|^2 - |\x_2 - 0.5|^2,
  \label{eq:solution_canuto}
\end{equation}
is the exact solution. We consider two different (although similar) \acp{QoI}, the average of the solution over the whole domain $\omega_1=\D(R)$ or over a (deterministic) subregion $ \omega_2=\widehat{D} \subset \D$ taken as $\widehat{D}=[0.5-\delta,0.5+\delta]^2$ with $\delta=0.125$, that is,
\begin{equation}
  Q_i(u):=\frac{1}{|\omega_i|}\int_{\omega_i}u
  \label{eq:qois}
\end{equation}
where $|\omega_i|$ denotes the area of $\omega_i$. Using \Eq{eq:solution_canuto} we get $Q_1=0.5 R^2$ and $Q_2=R^2 - 2\delta^2/3$. The \ac{PDF} of the $\TN$ can be written in terms of the error function and the variance can be computed easily using, e.g., \cite{truncatednormal}, which gives 
$\E(Q_1)=0.04531216540324139$ and $\E(Q_2)=0.08020766413981611$.

Because the \ac{MLMC} estimation given by \Eq{eq:mlmc_average} is a random variable so is the squared error $(\widetilde{Q}_L-\E(Q))^2$. Several realizations of the experiment will produce different results and its expectation error is determined by \eqref{eq:mlmc_error}. Therefore, to evaluate the error decay in \eqref{eq:error_decay} we approximate it by a sample average, running each experiment $K$ times, i.e., we compute the
error as
\begin{equation}
  \e_L^2 = \frac{1}{K} \sum_{k=1}^K (\widetilde{Q}^i_L-\E(Q))^2, 
  \label{eq:average_error}
\end{equation}
where $\widetilde{Q}^k_L$ is the result of the $k$-th realization of the experiment. The numerical experiments have been performed with $K=100$, whereas a value of $K=30$ was found to be sufficient in the context of hyperbolic conservation laws \cite{Mishra2012b}. Likewise we compute the averaged variances of the \ac{QoI} differences between levels as
\begin{equation}
 \mathcal{V}_l = \frac{1}{K} \sum_{k=1}^K V(Y_l^k)
  \label{eq:average_variance}
\end{equation}
where $V(Y_l^k)$ is given by \Eq{eq:sample_variance} for the $k$-th realization of the experiment. In practice, each realized \ac{MLMC} computation, here indexed by $k$, is computed with a different seed of the \ac{PRNG}. In the studies here we utilize the xoshiro256** algorithm \cite{PRNGxoshiro256} for generating pseudo random sequences, that has a 256 bit state. The algorithm is suited for this work as it is quick to advance the state, and statistically robust when seeded well. The detailed seeding procedure is not especially relevant to the results presented here, and those details are a bit beyond the scope of this paper, but the basic process can be given. The seeding process uses the recommended splitmix64 \ac{PRNG} \cite{PRNGsplitmix} to advance from an initial 64 bit seed. The outputs from this are used, together with physically generated random bits, to seed the state space of a Mersenne Twister \cite{Matsumoto1998}, specifically MT19937-64, and the outputs of the Mersenne Twister are used, together with physically generated random bits, to seed the initial state of the xoshiro256**. The state space of the generator has a period of $2^{256}-1$, meaning that after $2^{256}$ evaluations the state of the generator repeats, and there are no repetitions before that. The algorithm is adapted to parallelization by having easily computed a jump of $2^{128}$ states. Using this, jumps are performed to set the initial state for each sample. As a result, each process has an effective $2^{128}$ calls to the random number generator before reaching a state that was realized by another process, allowing for effective and efficient large scale computations. The jumps are done consistently, so that the initial state depends only on the seed and sample number, to allow reproducibility of the stochastic inputs regardless of the scheduling of order of computations.  

We consider $L=5$ and the number of samples per level are defined according to \Eq{eq:num_samples_ratio} with $\Gamma=7/2$. Each sample is computed taking $\artdom=[0,1]^2$ and $\T_0$ a Cartesian mesh of $8 \times 8 $ \acp{FE} and $\T_{\lev+1}$ is generated by dividing each element of $\T_\lev$ into four subelemements (two per direction). \Fig{fig:canuto_avg_err} shows the average error $\e_L$ as a function of the level for $N_5=\{3,6\}$. Even these low values of $N_5$, shows a good convergence of the level estimates, in agreement with the expected decay rate in \Ass{as:1}.

\begin{figure}
\centering
\includegraphics[scale=0.024]{QoIAveragesByLevelCombinedCanuto2d.pdf}
\caption{Averaged error $\e_L$ (given by \Eq{eq:average_error}) for each level and sample size (as from \Eq{eq:num_samples_ratio}) evaluated with $Q_1$ (Full Domain) and $Q_2$ (Partial Domain) from \Eq{eq:qois} compared to the theoretical expeted decay rate.}
\label{fig:canuto_avg_err}
\end{figure}

\begin{figure}
\centering
\includegraphics[scale=0.024]{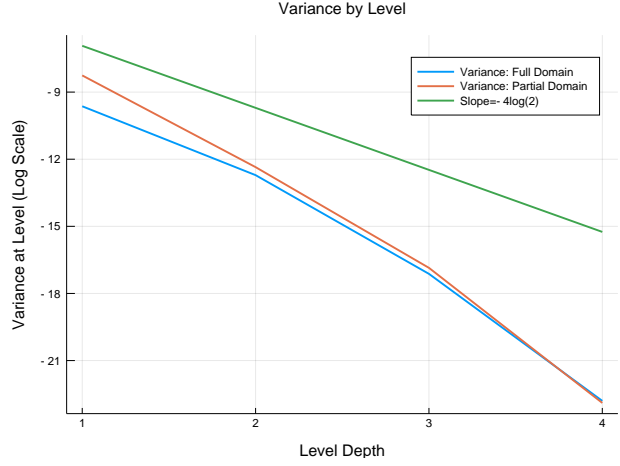}
\caption{Averaged sample variance between levels (given by \Eq{eq:average_variance}) evaluated with $Q_1$ (Full Domain) and $Q_2$ (Partial Domain) from \Eq{eq:qois}.}
\label{fig:canuto_var}
\end{figure}

\Fig{fig:canuto_var} shows the sample variance estimate \eqref{eq:average_variance} as a function of the level. The variance of the \ac{QoI} differences between levels decays with the level at the rate slightly higher than the one predicted by the theory in \cite{Teckentrup2013a}, verifying \Ass{as:3}. The determination of these variances is important in adaptive \ac{MLMC} methods as they are required to obtain the optimal number of samples per level in \Eq{eq:mlmc_num_samples}.

\begin{figure}
\centering
\includegraphics[scale=0.024]{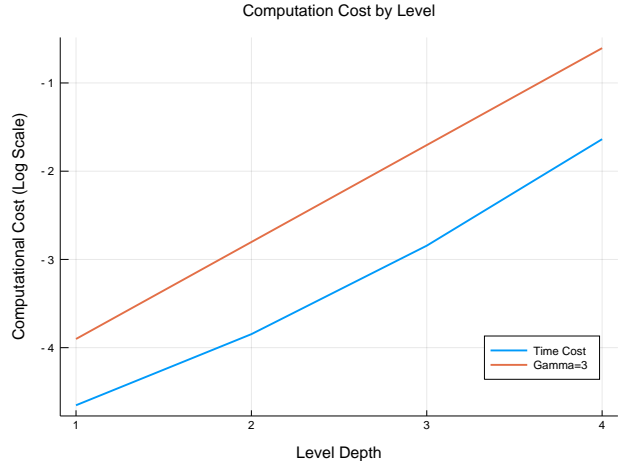}
\caption{Computational cost as a function of the level.}
\label{fig:canuto_cost}
\end{figure}

The number of samples in \Eq{eq:mlmc_num_samples} also depends on the complexity of sampling at each level  $C_\lev$ which is shown in \Fig{fig:canuto_cost}. This complexity is measured as the CPU time per sample averaged at each level, i.e., if $t^i_\lev$ is the elapsed CPU time in sample $i$ of level $\lev$, we estimate the complexity of sampling at each level as
\begin{equation}
  C_\lev=\frac{1}{N_\lev}\sum_{i=1}^{N_\lev} t^i_\lev.
\end{equation}
As it can be seen in \Fig{fig:canuto_cost}, the scaling is in agreement with \Ass{as:2}.

\begin{figure}
\centering
\includegraphics[scale=0.024]{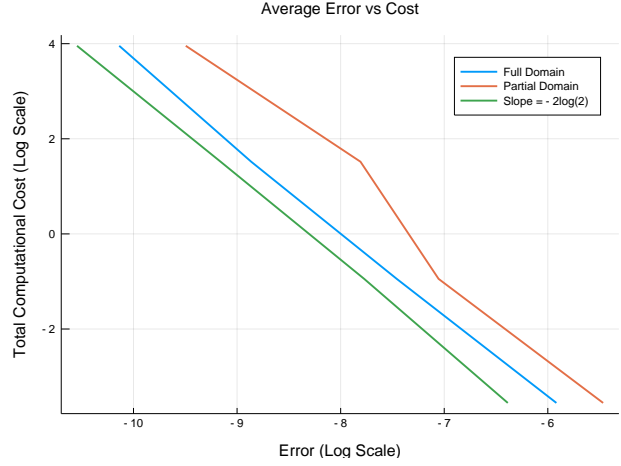}
\caption{Trade off between error and cost, evaluated with $Q_1$ (Full Domain) and $Q_2$ (Partial Domain) from \Eq{eq:qois}.}
\label{fig:canuto_samples}
\end{figure}

After verifying that the EMLMC algorithm shows an error and variance decay and cost increase that is in agreement with \Ass{as:1} to \Ass{as:3} and because these assumptions permit to obtain the complexity bound \Eq{eq:mlmc_complexity_bound} we expect the actual error to decrease with cost increase accordingly. The actual error decrease with respect to the total complexity, measured as $$C = \sum_{\lev=0}^{L}\sum_{i=1}^{N_\lev} t^i_\lev$$ is shown in \Fig{fig:canuto_samples} and the agreement with the scaling \Eq{eq:mlmc_complexity_bound} is evident. The important gain of \ac{MLMC} with respect to \ac{MC} in deterministic domains with random coefficients is also observed for the \ac{EMLMC} with this random geometry.

\subsection{Dealing with complex random geometries}
\label{subsec:complex}

The aim of this example is to show the robustness of the proposed approach to perform \ac{UQ} in complex random geometries. To this end, we consider a stochastic extension of a benchmark commonly used to test numerical methods designed to deal with interfaces described by a level-set \cite{hautefeuille_robust_2012,burman_cutfem:_2015}, the Poisson problem in a popcorn domain, a ball with spherical protuberances attached to its surface. Therefore, we aim to solve problem \Eq{eq:model} with $k=1$ and $f=-\nabla^2 u_{\rm ex}$ and Dirichlet boundary conditions $u_0=u_{\rm ex}$, where
\begin{equation}
u_{\rm ex} = \sin(k\|\x-\x_c\|)
\end{equation}
is defined in $\artdom=[0,1]^d$ with $\x_c=(1/2,1/2,1/2)$ and each sample is computed in a domain $\D(\w)$ defined through a level set function, i.e. $\x \in \D(\w)$ iff $\phi(\x,\w)<0$. The classical definition of the popcorn domain is
\begin{equation}
  \phi(\x) = \|\x - \x_0 \| - \rho_0 - \sum_{i=1}^n A \exp{ \|\x- \x_i\| / \sigma}
\end{equation}
with given parameters $\rho_0$ (the radius of the ball), $\x_i$ (the location of spherical protuberances), $A$ and $\sigma$ (which control the size of the protuberances) (see e.g. \cite{burman_cutfem:_2015}).

We consider a stochastic variation of this geometry defined by a random level set for $d=2,3$, although in this section we are primarily interested in the case $d=3$. This function is build to represent a random central ellipse (labelled as $0$) with smaller random ellipses attached to its boundary (labelled by $j$ with $1\leq j \leq n$) and is defined by
\begin{equation}
  \phi(\x,\w) =  \min_{j\in \{0:n\}}\left( \|\x-\x_j\|_{\bm{AR}(j)} - \rho_j\right),
  \label{eq:level_set_popcorn}
\end{equation}
where $x_j$ are random centers and $\rho_j$ random radius. The deformation from circles to ellipses is achieved by the norm
\begin{equation}
\|\y\|_{\bm{AR}(j)} = \sqrt{\left<(\bm{A}_j\bm{R}_j)^{-1}\y,(\bm{A}_j\bm{R}_j)^{-1}\y\right>}.
\end{equation}
which depends on streatching ($\bm{A}_j$) and rotation ($\bm{R}_j$) matrices defined as follows. Distortion matrices $\bm{A}_j$ are square, diagonal matrices with diagonal entries ($(\bm{A}_j)_{ii}$) defined by parameters $A^{\prime}_j(i)$ ($1\leq i \leq d$) 
\begin{equation}
  \label{eq:axis_normalization_def}
  (\bm{A}_j)_{ii} := \frac{A^{\prime}_j(i))}{\sqrt{\sum_{i=1}^d (A^{\prime}_j(i))^2}}. 
\end{equation}
Rotations are defined in each $(x_i,x_{i+1})$ plane ($1\leq i\leq d-1$)
through the rotation matrix $\bm{R}^{\prime}(\Theta_j(i))$ defined from random parameters $\Theta_j(i)$ as
\begin{equation}
  \label{eq:rotation_def1}
  \bm{R}^{\prime}(\Theta_j(i)) := \left(\begin{array}{cc}
   \cos(\Theta_j(i)) & -\sin(\Theta_j(i))\\
   \sin(\Theta_j(i)) & \cos(\Theta_j(i))\\
  \end{array}\right),
\end{equation}
Let $\bm{R}(\Theta_j(i))$, be the embedding of $\bm{R}^\prime(\Theta_j(i))$ inside $\R^d$, acting as the identity on the other $d-2$ coordinates of the vector. We then define the application of all rotations in $\R^d$ by one matrix as $\bm{R}_j$ defined by applying these rotations $\bm{R}(\Theta_j(1)),\cdots, \bm{R}(\Theta_j(d-1))$, in order, or as one single matrix denoted as,
\begin{equation}
  \label{eq:rotation_def3}
  \bm{R}_j := \prod_{i=1}^{d-1}\bm{R}(\Theta_j(i)) = \bm{R}(\Theta_j(d-1))\cdots\bm{R}(\Theta_j(1)) .
\end{equation}

\begin{figure}[!b]
  \centering
\includegraphics[width=0.45\textwidth]{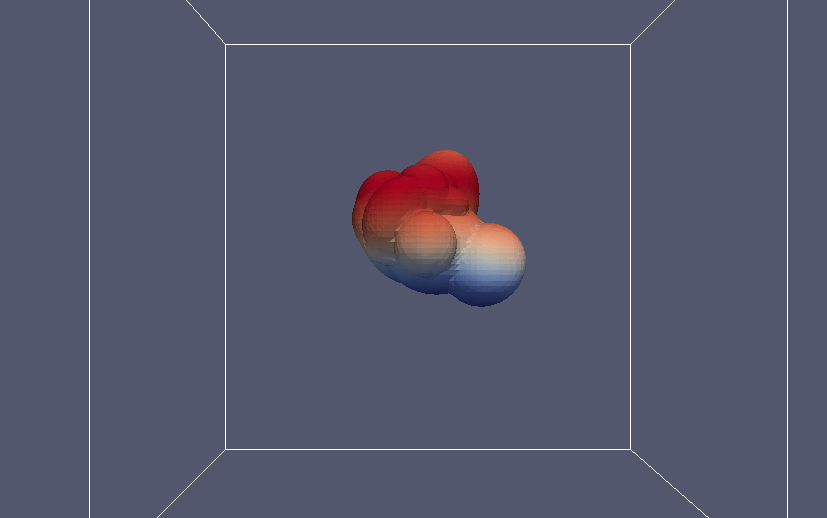} 
\includegraphics[width=0.45\textwidth]{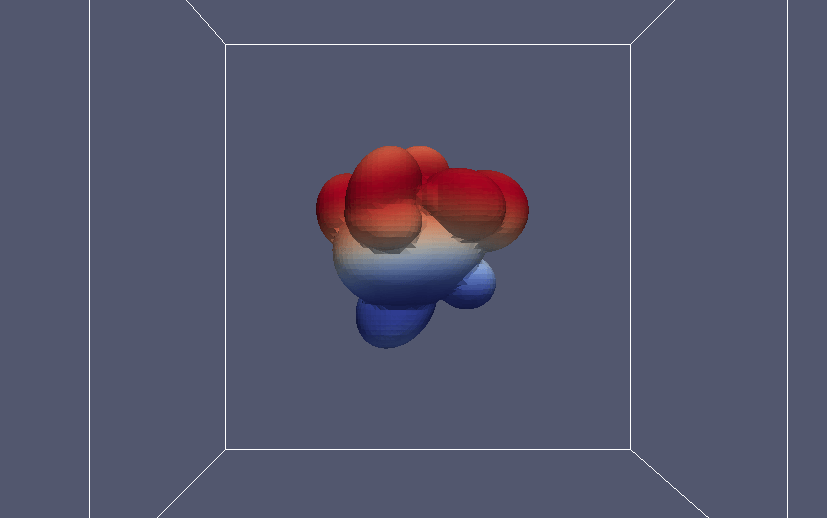} \\
\includegraphics[width=0.45\textwidth]{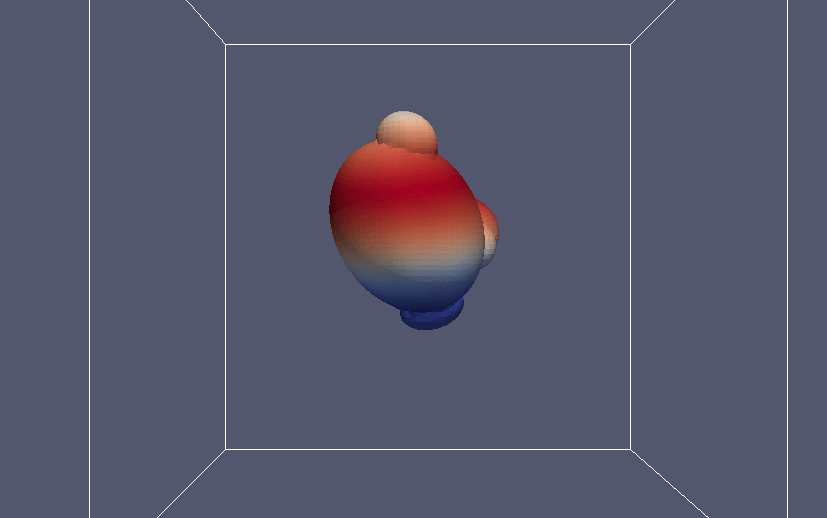}
\includegraphics[width=0.45\textwidth]{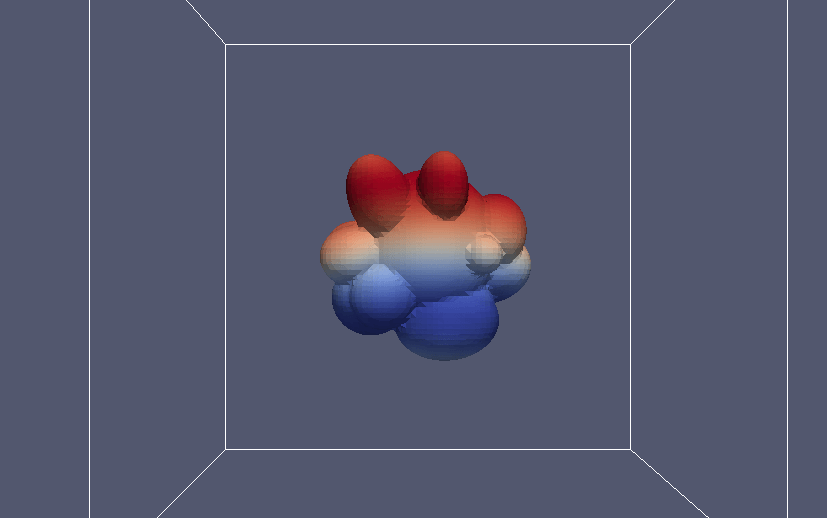} \\
\includegraphics[width=0.45\textwidth]{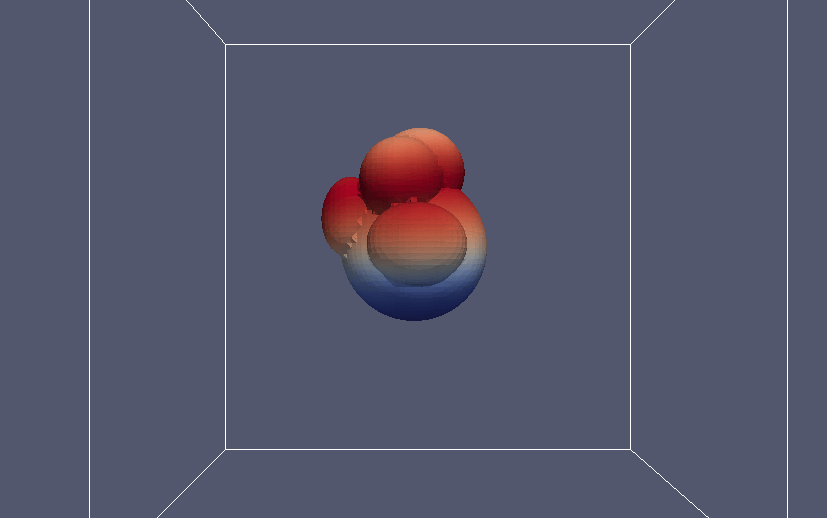} 
\includegraphics[width=0.45\textwidth]{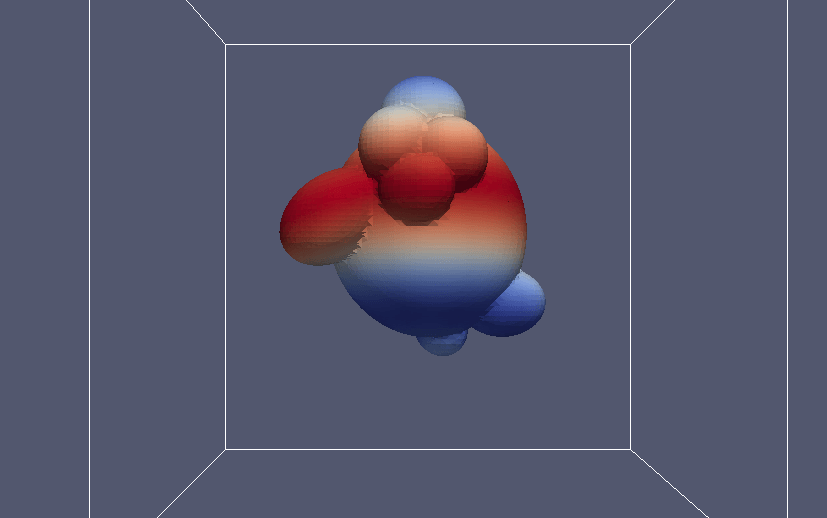} \\
  \caption{Some realizations of the solution of \Eq{eq:model} in a stochastic popcorn domain.}
  \label{fig:popcorn_samples}
\end{figure}

The random parameters describing this geometry are distributed according to
\begin{itemize}
\item $n \in \P(11)$
\item $A^{\prime}_j \in \U(0.8,1.3)^d$
\item $\x_0 \in \U(0.4,0.6)^d$
\item $\rho_0 \in \U(0.1,0.2)$
\item $\bm\Theta_j \in \U(0,2\pi)^{d-1}$
\item $\y_j \in \S^d$
\item $\x_j=\x_0+ \rho_0\bm{A}_0\bm{R}_0\y_j$
\item $\rho_j \in \U(0.03,0.1)$
\end{itemize}
where $\P$ and $\U$ denote a Poisson and uniform distributions, respectively, and $\S^d$ is a uniform distribution on the surface of a $d$-dimensional unit sphere. Note that using a level-set description permits to define complex domain through Boolean operations. In particular, given two level-set functions $\phi_1$ and $\phi_2$ describing $\Omega_1$ and $\Omega_2$ respectively, the function $\min(\phi_1,\phi_2)$ describes  $\Omega_1 \cup \Omega_2$.
A few samples of these stochastic domains are shown in \Fig{fig:popcorn_samples}, where the important variation of the morphology and size of the domains is apparent.

\begin{figure}[!b]
  \centering
  \includegraphics[width=0.7\textwidth]{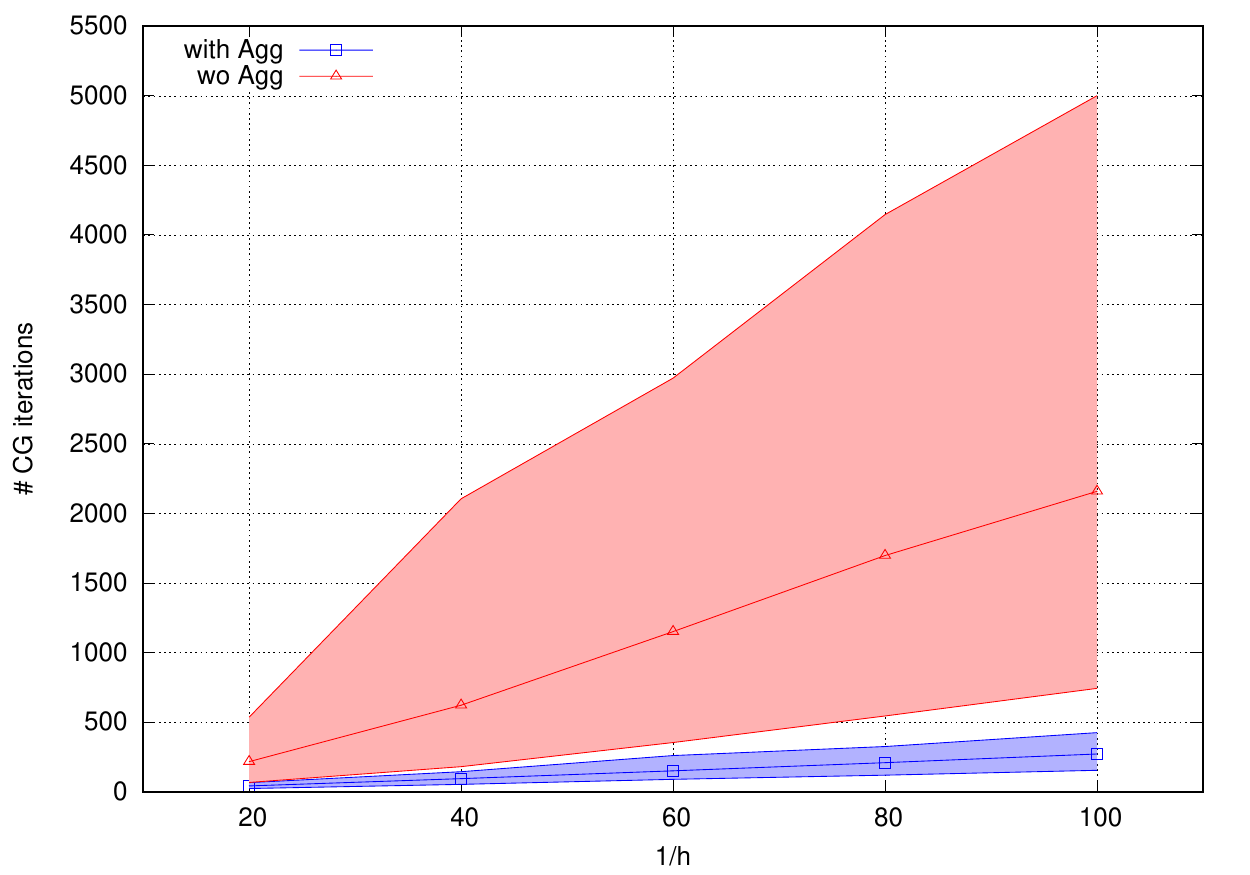}
  \caption{Number of \ac{CG} iterations required to converge to a tolerance of $10^{-8}$ in the relative norm the solution of each (2D) sample of the \ac{EMLMC} as a function of the mesh size on each level, with and without aggregation.}
  \label{fig:num_iter_2D}
\end{figure}

\begin{figure}[!b]
  \centering
  \includegraphics[width=0.45\textwidth]{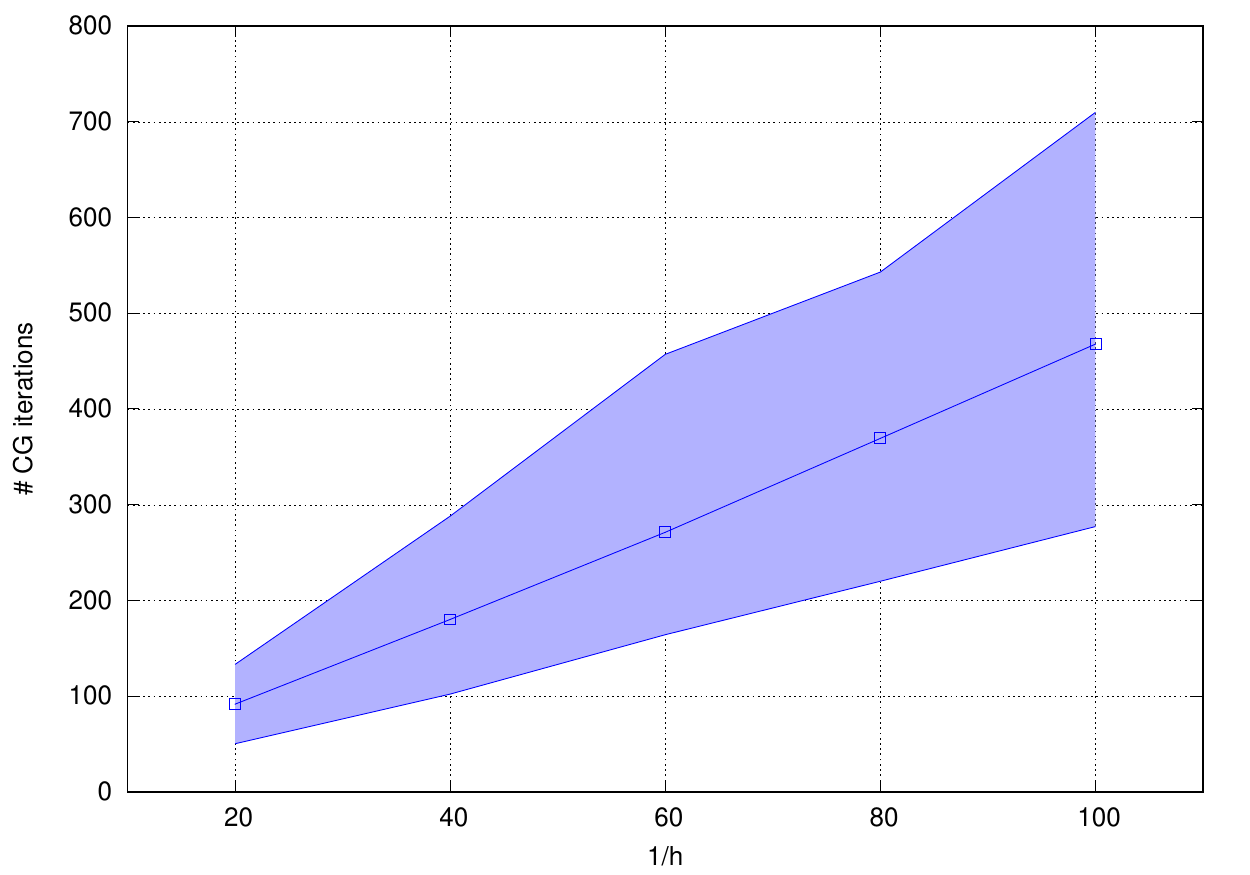}
  \includegraphics[width=0.45\textwidth]{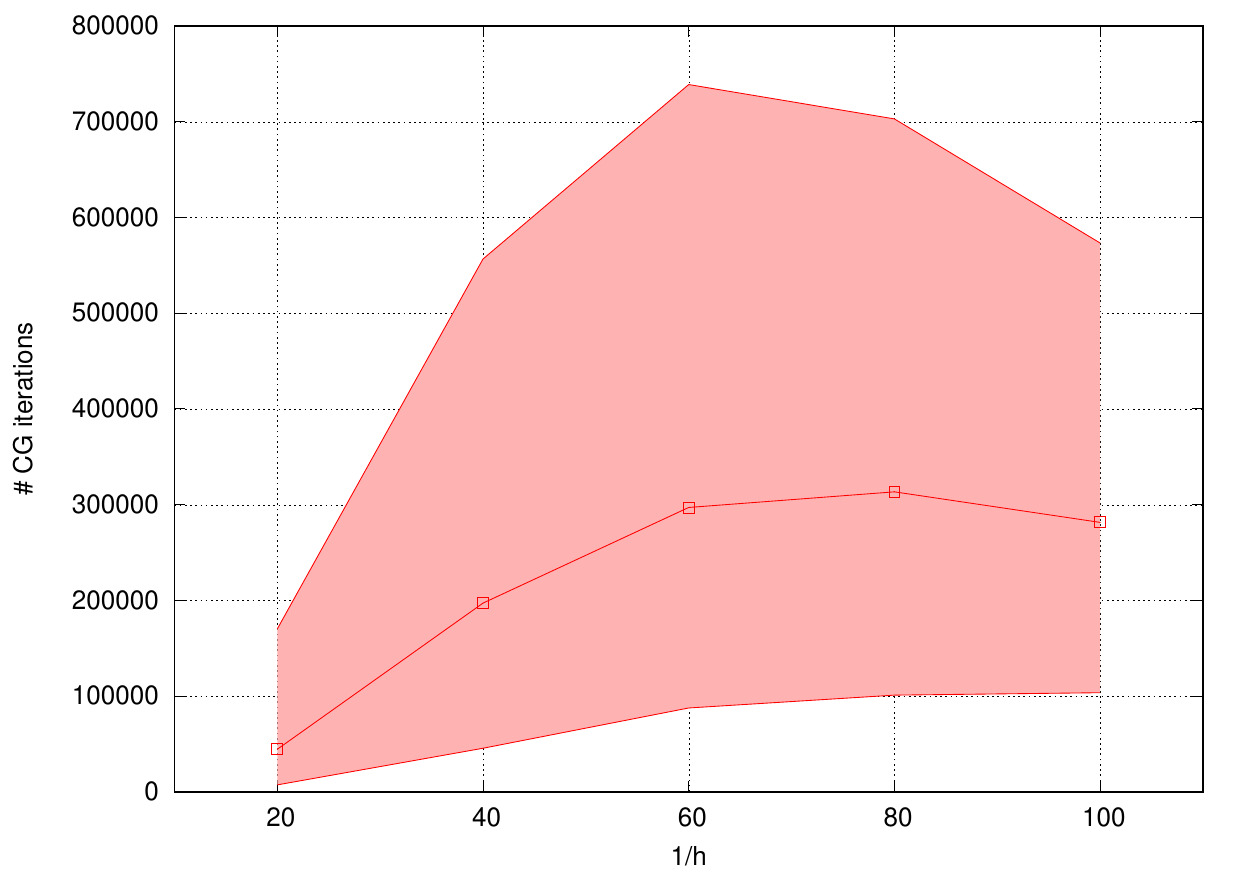}
  \caption{Number of \ac{CG} iterations required to converge to a tolerance of $10^{-8}$ in the relative norm the solution of each (3D) sample of the \ac{EMLMC} as a function of the mesh size on each level, (a) with and (b) without aggregation (note the difference in the scales).}
  \label{fig:num_iter_3D}
\end{figure}

An important problem with sampling methods, as described in \Sec{sec:intro} is the robustness of the solver employed to compute samples and this is an important feature of the \ac{agfem} employed herein and described in \Sec{sec:agg}. We performed a quantitative evaluation of this feature by running the \ac{EMLMC} method with and without the aggregation procedure using an iterative solver for the linear system arising from the \ac{FE} discretization. Due to the symmetry of problem \eqref{eq:model}, we use the \ac{CG} algorithm; as it is well-known, the number of iterations required to converge is proportional to the square root of the condition number of the system matrix. We compute statistics of the number of iterations (maximum, minimum and average) taking 1000 samples per level. The results obtained in a 2D version of the problem are shown in \Fig{fig:num_iter_2D} whereas those obtained in the 3D setting are shown in \Fig{fig:num_iter_3D}.

\subsection{Dealing with stochastic topologies}

In this section we consider a problem in  a stochastic domain where the randomness may induce changes in the topology. It consists of a square plate with two circular holes of uncertain position which may overlap to form a single one. We consider problem \eqref{eq:model} with Dirichlet boundary conditions on the left and right sides and zero Neumann boundary on the top and bottom sides and on the boundary of the internal holes. Physically, this problem corresponds to the heat transfer in a plate where the left (low) and right (high) temperatures are prescribed and heat flows through the plate at a rate that depends on the location and size of the holes. The natural quantity of interest is therefore 
\begin{equation}
  Q(u) :=\int_{x_1=0} \partial_{x_1} u,
  \label{eq:qoi_def_heat}
\end{equation}

The stochastic domain is defined, again, through a level set function $\phi(\x,\w) = \phi_1(\x,\w)+ \phi_2(\x,\w)$ where
\begin{equation}
  \phi_i(\x,\w) = \|\x - \x_i(\w) \| - \rho_i
\end{equation}
for $i=1,2$ represent the two interior circular holes (which do not conduct the heat) are distributed one along the upper half of the plate and one along the lower half of the plate. We study the influence of the radius on the average total heat flux throuogh the plate for random position of the holes. We conside three cases $\rho_i=0.18$, $\rho_i=0.2$ or $\rho_i=0.22$ for $i=1,2$. The first components of $\x_i$ for $i=1,2$ (holes' centers) are drawn uniformly from $[0.23,0.77]$. The upper circle has the $x_2$ coordinate of its center drawn uniformly from $[0.70,0.76]$, and the lower circle has the $x_2$ coordinate of its center drawn uniformly from $[0.24,0.30]$. The union of these circles form the non-conductive interior domain.

\begin{figure}[!b]
  \centering
  \includegraphics[width=0.45\textwidth]{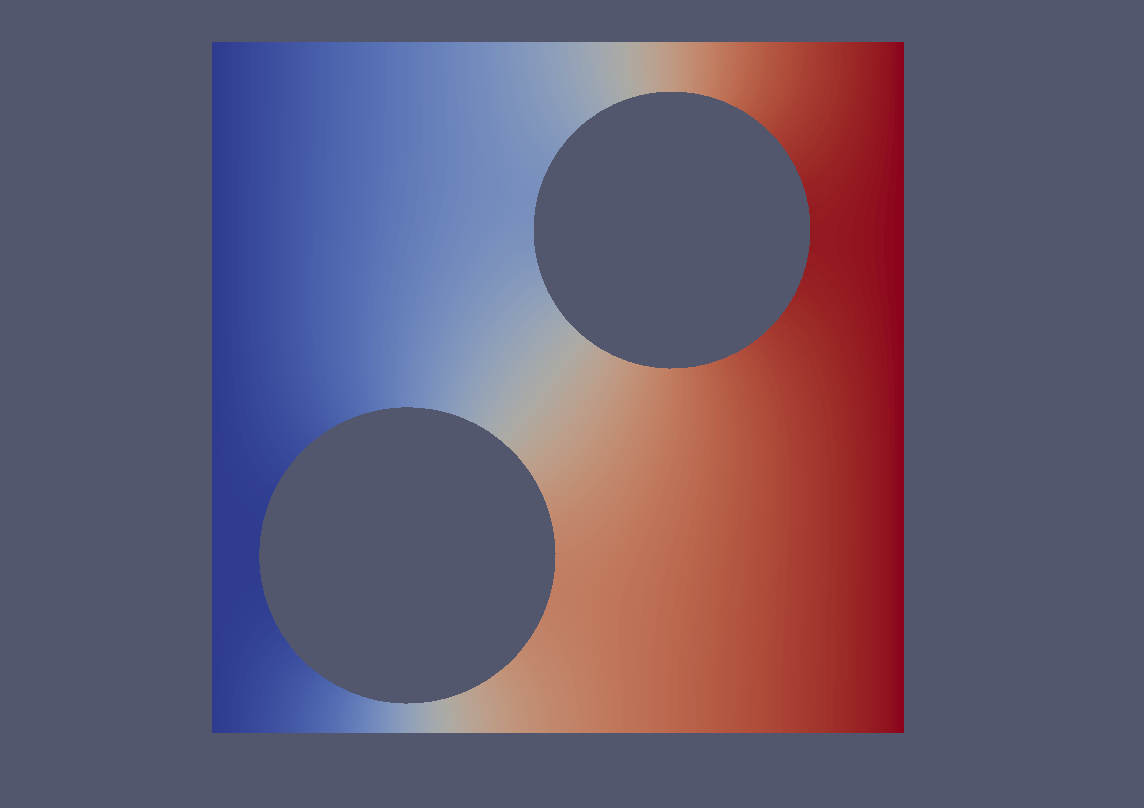} 
  \includegraphics[width=0.45\textwidth]{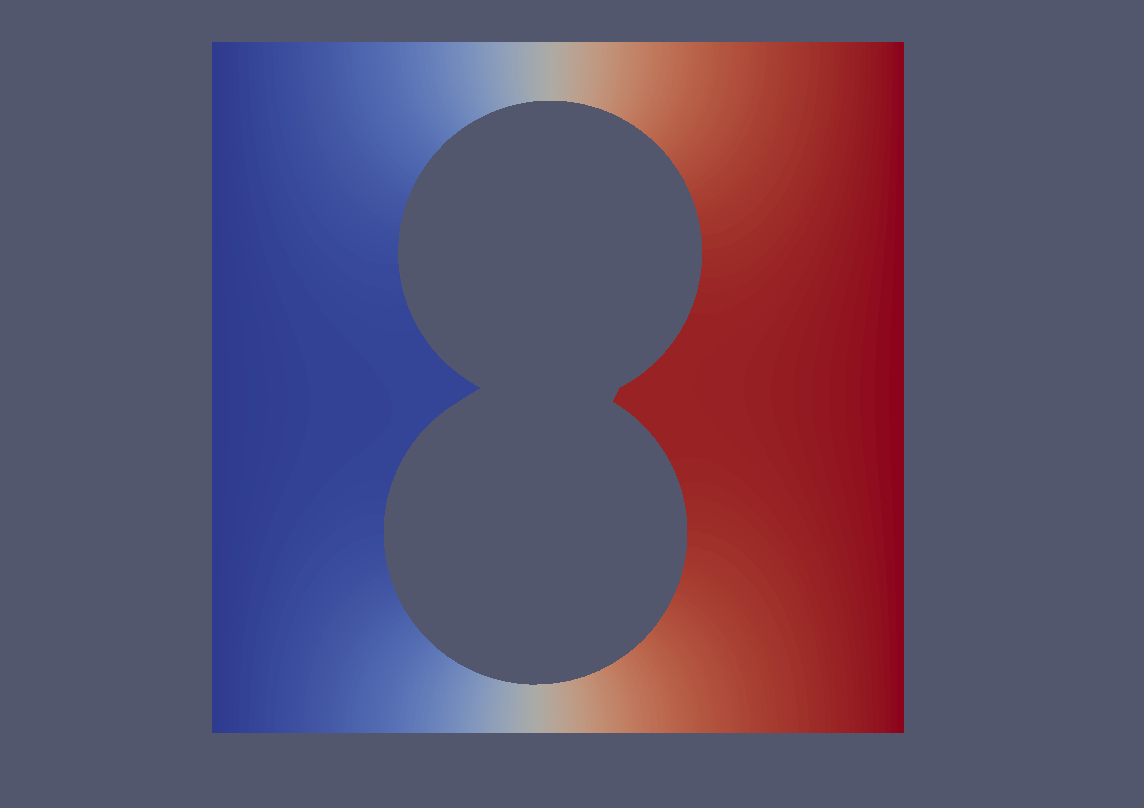} \\
  \includegraphics[width=0.45\textwidth]{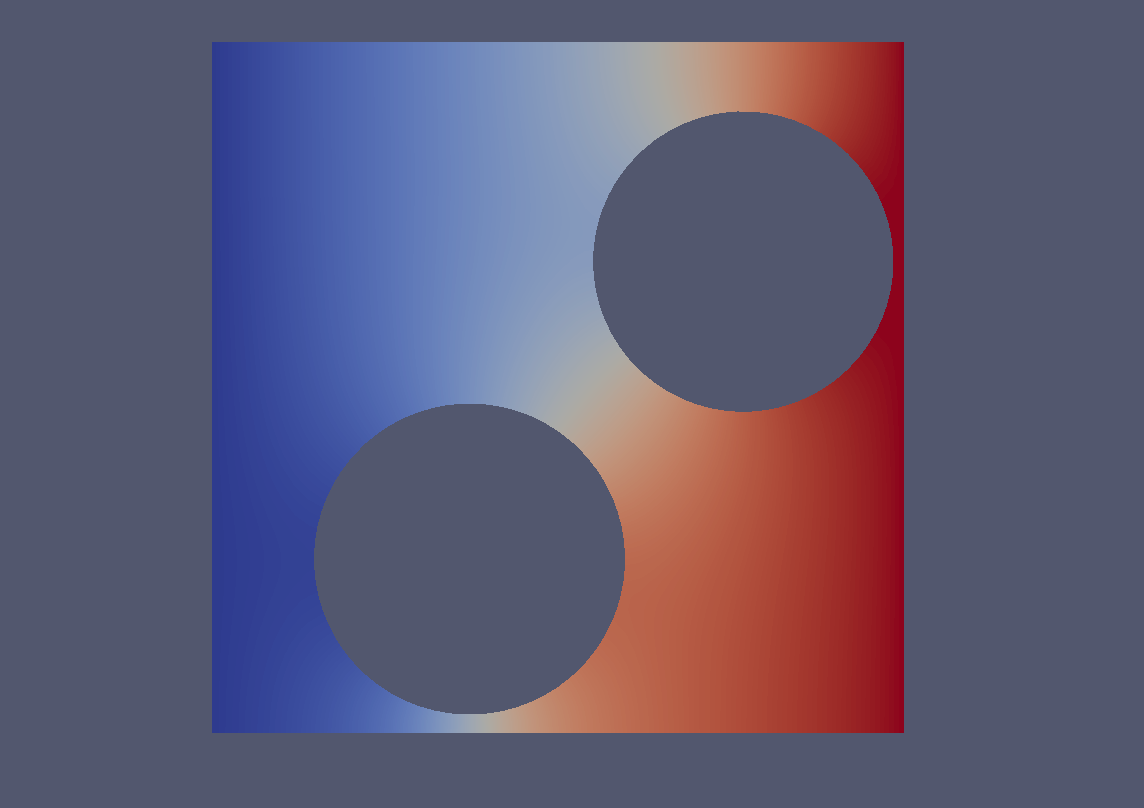}
  \includegraphics[width=0.45\textwidth]{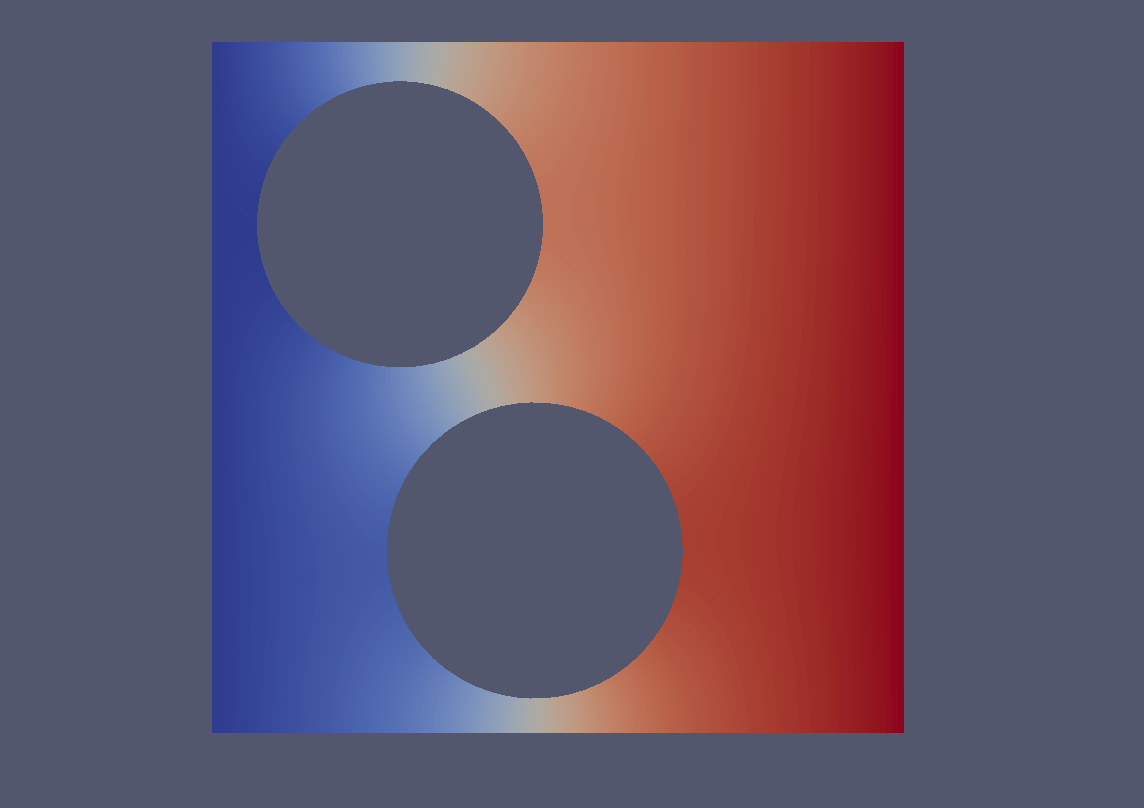} \\
  \includegraphics[width=0.45\textwidth]{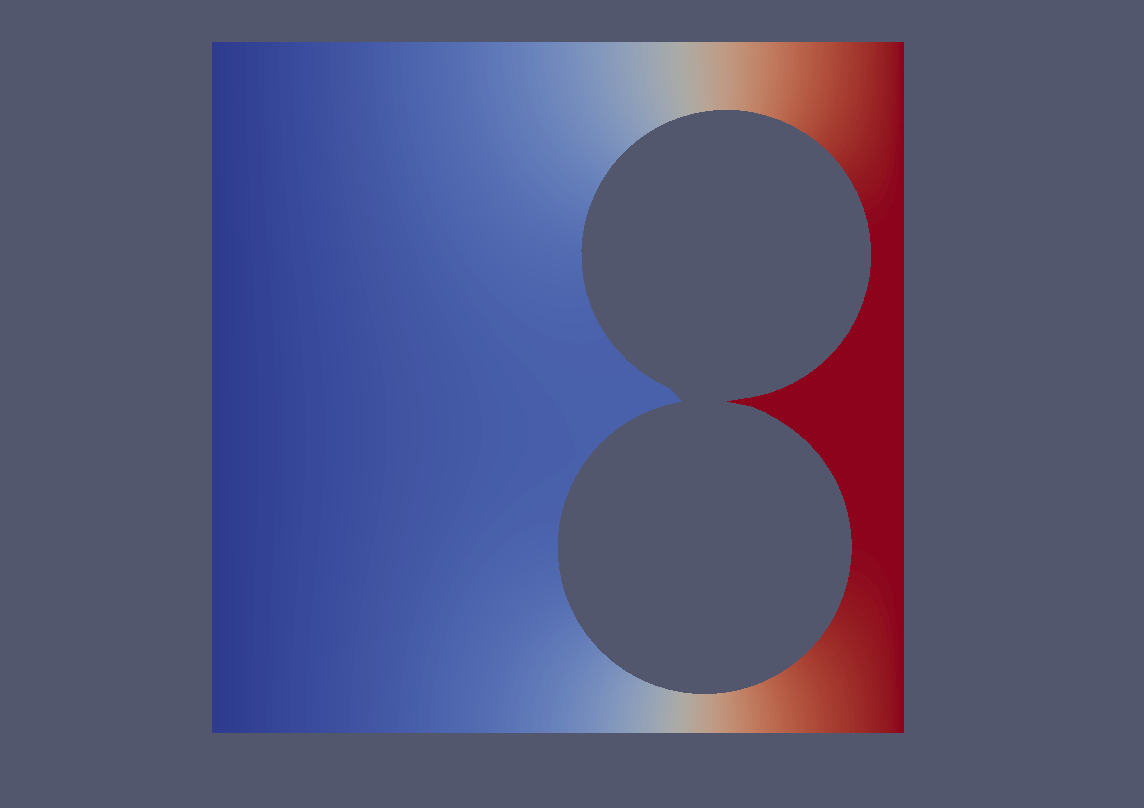} 
  \includegraphics[width=0.45\textwidth]{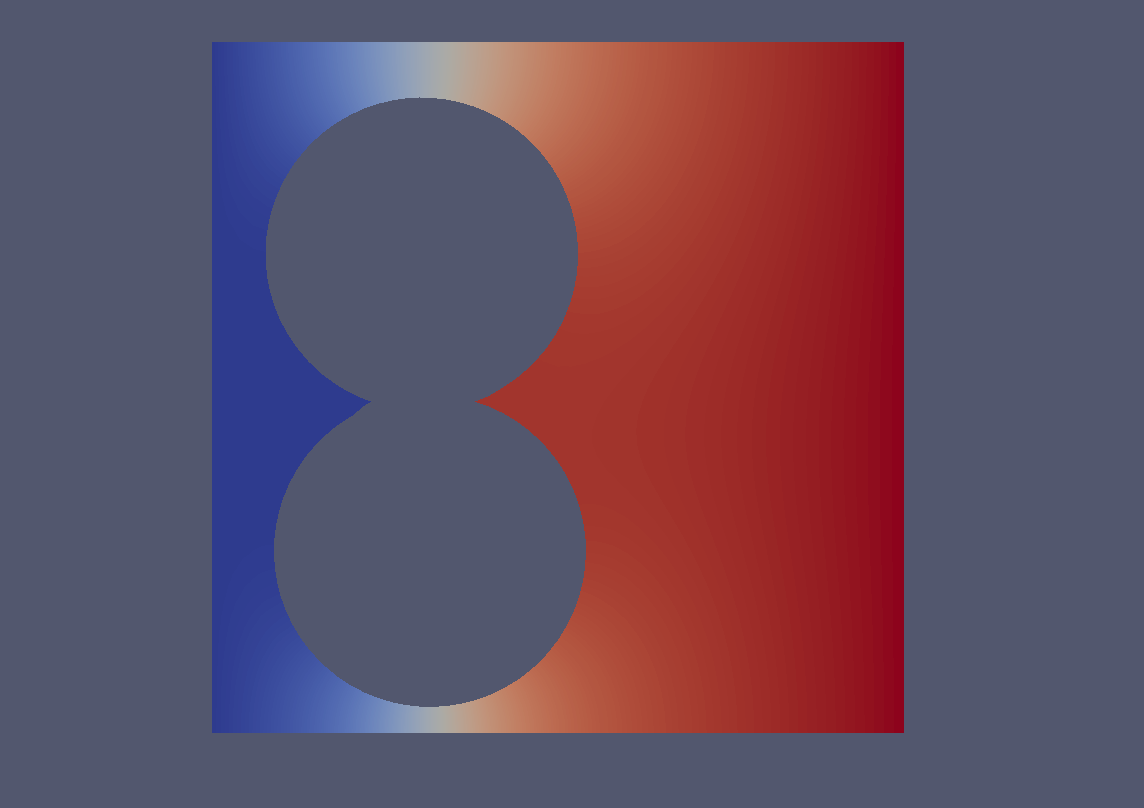} \\
  \caption{Some realizations of the solution of \Eq{eq:model} in a plate with two (stochastic) holes.}
  \label{fig:two_holes_geometry}
\end{figure}

We consider $L=5$ and the number of samples per level are defined according to \Eq{eq:num_samples_ratio} with $\Gamma=7/2$. Each sample is computed by taking $\artdom=[0,1]^2$ and $\T_0$ a Cartesian mesh of $8 \times 8 $ \acp{FE} and $\T_{\lev+1}$ is generated by dividing each element of $\T_\lev$ into four subelemements (two per direction). A few samples of the solution are shown in \Fig{fig:two_holes_geometry}, where the topology changes in the domain can be observed. The proposed \ac{EMLMC} method is able to automatically compute realizations using the same procedure followed in other cases, no special treatment of topology changes is required.

\begin{figure}[!b]
  \centering
  \includegraphics[scale=0.024]{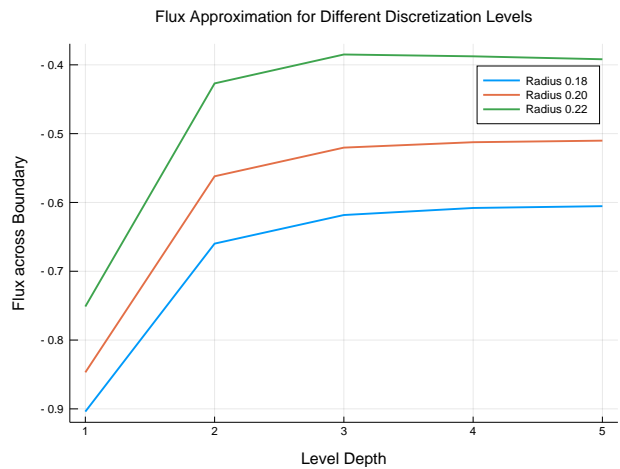}
  \caption{Average $Q(u)$ value as from \Eq{eq:qoi_def_heat} for three different radii at different levels of computation.}
  \label{fig:two_holes_flux}
\end{figure}

\Fig{fig:two_holes_flux} shows $\widetilde{Q}_L $ for the three different radii. As it can be seen, and is expected, when the size of the holes is increased, the effective area available for heat conduction is decreased, thus resulting in a reduction of the heat flux through the plate.

\section{Conclusions}\label{sec:conclusions}.

In this article we propose the \ac{EMLMC} method which consists of drawing samples in \ac{MLMC} using a robust embedd method, the \ac{agfem}, for the discretization of the problem. Once again, we emphasize that applying the standard \ac{MLMC} on complex geometries can be hard if body-fitted methods are used for the problem discretization due to the difficulty of the mesh-hierarchy generation, which tipically require human intervenation. On top of that, when dealing with random geometries the construction of the mesh-hierarchy might even be impossible using body-fitted techniques or may simply fail for a particular sample. Using embedded discretization methods on a mesh-hierarchy generated on a bounding-box circumvents these problems.

The resulting algorithm is a powerful method to perform \ac{UQ} on complex random domains thanks to the robustness provided by the \ac{agfem} algorithm. The numerical examples presented herein show that the error and variance decays and the complexity are similar to those observed in the standard \ac{MLMC} (with deterministic geometries and body-fitted discretizations) and are in line with theoretical expectations. Therefore, the \ac{EMLMC} shows the same cost reduction with respect to \ac{MC} than the standard \ac{MLMC} making it an excellent method to perform \ac{UQ} on complex random domains.

\bibliographystyle{plain}
\bibliography{arxiv}

\begin{thebibliography}{10}

\bibitem{truncatednormal}
{https://github.com/cossio/TruncatedNormal.jl}.

\bibitem{Babuska2001}
Ivo Babuska and Jan Chleboun.
\newblock {Effects of uncertainties in the Domain on the Solution of Neumann
  Boundary Value Problems}.
\newblock {\em Mathematics of Computation}, 71(240):1339--1370, 2001.

\bibitem{Babuska2003}
Ivo Babuska and Jan Chleboun.
\newblock {Effects of uncertainties in the domain on the solution of Dirichlet
  boundary value problems}.
\newblock {\em Numerische Mathematik}, 93(4):583--610, 2003.

\bibitem{Babuska2010}
Ivo Babuska, Fabio Nobile, and Ra\'{u}l Tempone.
\newblock {A Stochastic Collocation Method for Elliptic Partial Differential
  Equations with Random Input Data}.
\newblock {\em SIAM Review}, 52(3):317--355, 2010.

\bibitem{workflow}
Santiago Badia, Jerrad Hampton, and Javier Principe.
\newblock {A massively parallel implementation of Multilevel Monte Carlo
  methods}.
\newblock {\em In preparation}, 2019.

\bibitem{Badia2019}
Santiago Badia and Alberto~F. Mart{\'{i}}n.
\newblock {A tutorial-driven introduction to the parallel finite element
  library FEMPAR v1.0.0}.
\newblock {\em arXiv:1908.00891}, 2019.

\bibitem{fempar}
Santiago Badia, Alberto~F. Mart{\'{i}}n, and Javier Principe.
\newblock {FEMPAR: An Object-Oriented Parallel Finite Element Framework}.
\newblock {\em Archives of Computational Methods in Engineering},
  25(2):195--271, oct 2018.

\bibitem{Badia2018b}
Santiago Badia, Alberto~F. Martin, and Francesc Verdugo.
\newblock {Mixed Aggregated Finite Element Methods for the Unfitted
  Discretization of the Stokes Problem}.
\newblock {\em SIAM Journal on Scientific Computing}, 40(6):B1541--B1576, jan
  2018.

\bibitem{Badia2017a}
Santiago Badia, Francesc Verdugo, and Alberto~F. Mart{\'{i}}n.
\newblock {The aggregated unfitted finite element method for elliptic
  problems}.
\newblock {\em Computer Methods in Applied Mechanics and Engineering},
  336:533--553, jul 2018.

\bibitem{Barth2011}
Andrea Barth, Christoph Schwab, and Nathaniel Zollinger.
\newblock {Multi-level Monte Carlo Finite Element method for elliptic PDEs with
  stochastic coefficients}.
\newblock {\em Numerische Mathematik}, 119(1):123--161, sep 2011.

\bibitem{PRNGxoshiro256}
David Blackman and Sebastiano Vigna.
\newblock Scrambled linear pseudorandom number generators.
\newblock {\em arXiv preprint arXiv:1805.01407}, 2018.

\bibitem{Blyth2003}
M~G Blyth and C~Pozrikidis.
\newblock {Heat conduction across irregular and fractal-like surfaces}.
\newblock {\em International Journal of Heat and Mass Transfer},
  46(8):1329--1339, 2003.

\bibitem{burman_cutfem:_2015}
Erik Burman, Susanne Claus, P~Hansbo, M~G Larson, and Andr{\'{e}} Massing.
\newblock {{CutFEM}: {Discretizing} Geometry and Partial Differential
  Equations}.
\newblock {\em International Journal for Numerical Methods in Engineering},
  104(7):472--501, 2015.

\bibitem{Cajueiro1999}
D~O Cajueiro, V~A {de A Sampaio}, C~M~C de~Castilho, and R~F~S Andrade.
\newblock {Fractal properties of equipotentials close to a rough conducting
  surface}.
\newblock {\em Journal of Physics: Condensed Matter}, 11(26):4985, 1999.

\bibitem{Canuto2007}
Claudio Canuto and Tomas Kozubek.
\newblock {A fictitious domain approach to the numerical solution of PDEs in
  stochastic domains}.
\newblock {\em Numerische Mathematik}, 107(2):257--293, aug 2007.

\bibitem{Castrillon-Candas2016}
Julio~E Castrill{\'{o}}n-Cand{\'{a}}s, Fabio Nobile, and Ra{\'{u}}l~F Tempone.
\newblock {Analytic regularity and collocation approximation for elliptic PDEs
  with random domain deformations}.
\newblock {\em Computers and Mathematics with Applications}, 71:1173--1197,
  2016.

\bibitem{Charrier2013}
Julia Charrier, Robert Scheichl, and Aretha~Leonore Teckentrup.
\newblock {Finite element error analysis of elliptic PDEs with random
  coefficients and its application to multilevel Monte Carlo methods}.
\newblock {\em SIAM Journal on Numerical Analysis}, 51(1):322--352, 2013.

\bibitem{Chaudhry20181127}
Jehanzeb~H Chaudhry, Nathanial Burch, and Donald Estep.
\newblock {Efficient Distribution Estimation and Uncertainty Quantification for
  Elliptic Problems on Domains with Stochastic Boundaries}.
\newblock {\em SIAM/ASA Journal on Uncertainty Quantification},
  6(3):1127--1150, jan 2018.

\bibitem{Cliffe2011}
K~A Cliffe, M~B Giles, R~Scheichl, A~L Teckentrup, C~W Oosterlee, and A~K
  {Borzi A Cliffe}.
\newblock {Multilevel Monte Carlo methods and applications to elliptic PDEs
  with random coefficients}.
\newblock {\em Comput Visual Sci}, 14:3--15, 2011.

\bibitem{Collier2015}
Nathan Collier, Abdul-Lateef Haji-Ali, Fabio Nobile, Erik von Schwerin, and
  Ra{\'{u}}l Tempone.
\newblock {A continuation multilevel Monte Carlo algorithm}.
\newblock {\em BIT Numerical Mathematics}, 55(2):399--432, 2015.

\bibitem{Dambrine2016921}
M~Dambrine, I~Greff, H~Harbrecht, and B~Puig.
\newblock {Numerical Solution of the Poisson Equation on Domains with a Thin
  Layer of Random Thickness}.
\newblock {\em SIAM Journal on Numerical Analysis}, 54(2):921--941, jan 2016.

\bibitem{Dambrine2017943}
M~Dambrine, I~Greff, H~Harbrecht, and B~Puig.
\newblock {Numerical solution of the homogeneous Neumann boundary value problem
  on domains with a thin layer of random thickness}.
\newblock {\em Journal of Computational Physics}, 330:943--959, feb 2017.

\bibitem{Deffeyes1982}
K.~S. Deffeyes, B.~D. Ripley, and G.~S. Watson.
\newblock {Stochastic geometry in petroleum geology}.
\newblock {\em Journal of the International Association for Mathematical
  Geology}, 14(5):419--432, oct 1982.

\bibitem{Doi91}
Akio Doi and Akio Koide.
\newblock {An Efficient Method of Triangulating Equi-Valued Surfaces by Using
  Tetrahedral Cells}.
\newblock {\em IEICE Transactions on Information and Systems},
  (E74-D):214--224, 1991.

\bibitem{Fyrillas2001}
M~M Fyrillas and C~Pozrikidis.
\newblock {Conductive heat transport across rough surfaces and interfaces
  between two conforming media}.
\newblock {\em International Journal of Heat and Mass Transfer},
  44(9):1789--1801, 2001.

\bibitem{ghanem1991stochastic}
R~Ghanem and P~D Spanos.
\newblock {\em {Stochastic Finite Elements: A Spectral Approach}}.
\newblock Springer-Verlag New York, 1991.

\bibitem{Giles2008}
Michael~B. Giles.
\newblock {Multilevel Monte Carlo Path Simulation}.
\newblock {\em Operations Research}, 56(3):607--617, jun 2008.

\bibitem{Giles2012}
Michael~B. Giles and Christoph Reisinger.
\newblock {Stochastic Finite Differences and Multilevel Monte Carlo for a Class
  of SPDEs in Finance}.
\newblock {\em SIAM Journal on Financial Mathematics}, 3(1):572--592, jan 2012.

\bibitem{Gordon2014}
Andrew Gordon and Catherine~E Powell.
\newblock {A Preconditioner for Fictitious Domain Formulations of Elliptic PDEs
  on Uncertain Parameterized Domains}.
\newblock {\em SIAM/ASA Journal on Uncertainty Quantification}, 2(1):622--646,
  jan 2014.

\bibitem{Haji-Ali2016}
Abdul-Lateef Haji-Ali, Fabio Nobile, Erik {Von Schwerin}, and Ra{\'{u}}l
  Tempone.
\newblock {Optimization of mesh hierarchies in multilevel Monte Carlo
  samplers}.
\newblock 4:76--112, 2016.

\bibitem{hamrock2004fundamentals}
Bernard~J Hamrock, Steven~R Schmid, and Bo~O Jacobson.
\newblock {\em {Fundamentals of fluid film lubrication}}.
\newblock CRC press, 2004.

\bibitem{Harbrecht2016}
H~Harbrecht, M.~Peters, and M~Siebenmorgen.
\newblock {Analysis of the domain mapping method for elliptic diffusion
  problems on random domains}.
\newblock {\em Numerische Mathematik}, 134(4):823--856, 2016.

\bibitem{Harbrecht2013}
Helmut Harbrecht and Jingzhi Li.
\newblock {First order second moment analysis for stochastic interface problems
  based on low-rank approximation}.
\newblock {\em ESAIM: Mathematical Modelling and Numerical Analysis},
  47(5):1533--1552, sep 2013.

\bibitem{Harbrecht2008}
Helmut Harbrecht, Reinhold Schneider, and Christoph Schwab.
\newblock {Sparse second moment analysis for elliptic problems in stochastic
  domains}.
\newblock {\em Numerische Mathematik}, 109(3):385--414, 2008.

\bibitem{hautefeuille_robust_2012}
Martin Hautefeuille, Chandrasekhar Annavarapu, and John~E Dolbow.
\newblock {Robust imposition of {\{}Dirichlet{\}} boundary conditions on
  embedded surfaces}.
\newblock {\em International Journal for Numerical Methods in Engineering},
  90(1):40--64, apr 2012.

\bibitem{hutchings2017tribology}
Ian Hutchings and Philip Shipway.
\newblock {\em {Tribology: friction and wear of engineering materials}}.
\newblock Butterworth-Heinemann, 2017.

\bibitem{Kebaier2005}
Ahmed Kebaier.
\newblock {Statistical Romberg extrapolation: A new variance reduction method
  and applications to option pricing}.
\newblock {\em The Annals of Applied Probability}, 15(4):2681--2705, nov 2005.

\bibitem{Lang20131031}
Christapher Lang, Alireza Doostan, Kurt Maute, Christapher Lang, Alireza
  Doostan, and {\textperiodcentered}~K Maute.
\newblock {Extended stochastic FEM for diffusion problems with uncertain
  material interfaces}.
\newblock {\em Computational Mechanics}, 51(6):1031--1049, jun 2013.

\bibitem{lemaitre2010spectral}
Olivier~P. {Le Maitre} and Omar~M. Knio.
\newblock {\em {Spectral Methods for Uncertainty Quantification. With
  Applications to Computational Fluid Dynamics}}.
\newblock Springer, 2010.

\bibitem{Luo2008}
Xian Luo, Martin~R Maxey, and George~Em Karniadakis.
\newblock {Smoothed profile method for particulate flows: Error analysis and
  simulations}.
\newblock {\em Journal of Computational Physics}, 228(5):1750--1769, 2009.

\bibitem{Matsumoto1998}
Makoto Matsumoto and Takuji Nishimura.
\newblock {Mersenne twister: a 623-dimensionally equidistributed uniform
  pseudo-random number generator}.
\newblock {\em ACM Transactions on Modeling and Computer Simulation},
  8(1):3--30, 1998.

\bibitem{Mishra2012b}
S~Mishra and Ch. Schwab.
\newblock {Sparse tensor multi-level Monte Carlo finite volume methods for
  hyperbolic conservation laws with random initial data}.
\newblock {\em Mathematics of Computation}, 81(280):1979--2018, 2012.

\bibitem{Mishra2012a}
S.~Mishra, Ch. Schwab, and J.~Sukys.
\newblock {Multi-level Monte Carlo finite volume methods for nonlinear systems
  of conservation laws in multi-dimensions}.
\newblock {\em Journal of Computational Physics}, 231(8):3365--3388, apr 2012.

\bibitem{Mishra2012}
S.~Mishra, Ch. Schwab, J~Sukys, and J.~Sukys.
\newblock {Multilevel Monte Carlo Finite Volume Methods for Shallow Water
  Equations with Uncertain Topography in Multi-dimensions}.
\newblock {\em SIAM Journal on Scientific Computing}, 34(6):B761--B784, jan
  2012.

\bibitem{Moes1999}
Nicolas Mo{\"{e}}s, John Dolbow, and Ted Belytschko.
\newblock {A finite element method for crack growth without remeshing}.
\newblock {\em International Journal for Numerical Methods in Engineering},
  46(1):131--150, sep 1999.

\bibitem{Mohan2011874}
P~Surya Mohan, Prasanth~B Nair, and Andy~J Keane.
\newblock {Stochastic projection schemes for deterministic linear elliptic
  partial differential equations on random domains}.
\newblock {\em International Journal for Numerical Methods in Engineering},
  85(7):874--895, feb 2011.

\bibitem{Nosonovsky2008}
Michael Nosonovsky and Bharat Bhushan.
\newblock {Biologically inspired surfaces: Broadening the scope of roughness}.
\newblock {\em Advanced Functional Materials}, 18(6):843--855, mar 2008.

\bibitem{Nouy20113066}
A~Nouy, M~Chevreuil, and E~Safatly.
\newblock {Fictitious domain method and separated representations for the
  solution of boundary value problems on uncertain parameterized domains}.
\newblock {\em Computer Methods in Applied Mechanics and Engineering},
  200(45-46):3066--3082, oct 2011.

\bibitem{Nouy20101312}
A~Nouy and A~Cl{\'{e}}ment.
\newblock {eXtended Stochastic Finite Element Method for the numerical
  simulation of heterogeneous materials with random material interfaces}.
\newblock {\em International Journal for Numerical Methods in Engineering},
  83(10):1312--1344, sep 2010.

\bibitem{Nouy2008}
A~Nouy, A~Cl{\'{e}}ment, F~Schoefs, and N~Mo{\"{e}}s.
\newblock {An extended stochastic finite element method for solving stochastic
  partial differential equations on random domains}.
\newblock {\em Computer Methods in Applied Mechanics and Engineering},
  197(51-52):4663--4682, 2008.

\bibitem{Nouy2007}
Anthony Nouy, Franck Schoefs, and Nicolas Mo{\"{e}}s.
\newblock {X-SFEM, a computational technique based on X-FEM to deal with random
  shapes}.
\newblock {\em European Journal of Computational Mechanics}, 16(2):277--293,
  2007.

\bibitem{Park2012}
Sang~Woo Park, Marcos Intaglietta, and Daniel~M Tartakovsky.
\newblock {Impact of endothelium roughness on blood flow}.
\newblock {\em Journal of Theoretical Biology}, 300:152--160, 2012.

\bibitem{Parussini2010}
Lucia Parussini, Valentino Pediroda, and Carlo Poloni.
\newblock {Prediction of geometric uncertainty effects on Fluid Dynamics by
  Polynomial Chaos and Fictitious Domain method}.
\newblock {\em Computers and Fluids}, 39(1):137--151, 2010.

\bibitem{Patz2013}
Torben P{\"{a}}tz, Robert~M Kirby, and Tobias Preusser.
\newblock {Ambrosio-Tortorelli Segmentation of Stochastic Images: Model
  Extensions, Theoretical Investigations and Numerical Methods}.
\newblock {\em International Journal of Computer Vision}, 103(2):190--212, jun
  2013.

\bibitem{Pisaroni2017}
M.~Pisaroni, F.~Nobile, and P.~Leyland.
\newblock {A Continuation Multi Level Monte Carlo (C-MLMC) method for
  uncertainty quantification in compressible inviscid aerodynamics}.
\newblock {\em Computer Methods in Applied Mechanics and Engineering},
  326:20--50, 2017.

\bibitem{Richardson1971}
S.~Richardson.
\newblock {A model for the boundary condition of a porous material. Part 2}.
\newblock {\em Journal of Fluid Mechanics}, 49(2):327--336, 1971.

\bibitem{Sarkar1993}
Kausik Sarkar and Charles Meneveau.
\newblock {Gradients of potential fields on rough surfaces: Perturbative
  calculation of the singularity distribution function {\{}f($\backslash$){\}}
  for small surface dimension}.
\newblock {\em Physical Review E}, 47(2):957--966, 1993.

\bibitem{Schoefs2016}
Franck Schoefs, Mathilde Chevreuil, Olivier Pasqualini, and Mika{\"{e}}l
  Cazuguel.
\newblock {Partial safety factor calibration from stochastic finite element
  computation of welded joint with random geometries}.
\newblock 2016.

\bibitem{sethian1999level}
James~Albert Sethian.
\newblock {\em {Level set methods and fast marching methods: evolving
  interfaces in computational geometry, fluid mechanics, computer vision, and
  materials science}}, volume~3.
\newblock Cambridge university press, 1999.

\bibitem{PRNGsplitmix}
Guy~L Steele~Jr, Doug Lea, and Christine~H Flood.
\newblock Fast splittable pseudorandom number generators.
\newblock In {\em ACM SIGPLAN Notices}, volume~49, pages 453--472. ACM, 2014.

\bibitem{Stefanou2009}
G.~Stefanou, A.~Nouy, and A.~Clement.
\newblock {Identification of random shapes from images through polynomial chaos
  expansion of random level set functions}.
\newblock {\em International Journal for Numerical Methods in Engineering},
  79(2):127--155, jul 2009.

\bibitem{Stefanou2004}
G~Stefanou and M~Papadrakakis.
\newblock {Stochastic finite element analysis of shells with combined random
  material and geometric properties}.
\newblock {\em Computer Methods in Applied Mechanics and Engineering},
  193(1-2):139--160, jan 2004.

\bibitem{sukumar_modeling_2001}
N~Sukumar, D~L Chopp, N~Mo{\"{e}}s, and T~Belytschko.
\newblock {Modeling holes and inclusions by level sets in the extended
  finite-element method}.
\newblock {\em Computer Methods in Applied Mechanics and Engineering},
  190(46–47):6183--6200, 2001.

\bibitem{Tartakovsky2006}
Daniel~M Tartakovsky and Dongbin Xiu.
\newblock {Stochastic analysis of transport in tubes with rough walls}.
\newblock {\em Journal of Computational Physics}, 217(1):248--259, 2006.

\bibitem{Taylor1971}
G.~I. Taylor.
\newblock {A model for the boundary condition of a porous material. Part 1}.
\newblock {\em Journal of Fluid Mechanics}, 49(02):319, sep 1971.

\bibitem{Teckentrup2013a}
A.~L. Teckentrup, R.~Scheichl, M.~B. Giles, and E.~Ullmann.
\newblock {Further analysis of multilevel Monte Carlo methods for elliptic PDEs
  with random coefficients}.
\newblock {\em Numer. Math}, 125:569--600, 2013.

\bibitem{Verdugo2019}
Francesc Verdugo, Alberto~F. Mart{\'{i}}n, and Santiago Badia.
\newblock {Distributed-memory parallelization of the aggregated unfitted finite
  element method}.
\newblock {\em Computer Methods in Applied Mechanics and Engineering},
  357:112583, dec 2019.

\bibitem{Xiu2003}
Dongbin Xiu and George~Em Karniadakis.
\newblock {Modeling uncertainty in flow simulations via generalized polynomial
  chaos}.
\newblock {\em Journal of Computational Physics}, 187(1):137--167, 2003.

\bibitem{Xiu2006}
Dongbin Xiu and D.M. Tartakovsky.
\newblock {Numerical methods for differential equations in random domains}.
\newblock {\em SIAM Journal on Scientific Computing}, 28(3):1167--1185, jan
  2006.

\bibitem{Zayernouri2013}
Mohsen Zayernouri, Sang-Woo Park, Daniel~M Tartakovsky, and George~Em
  Karniadakis.
\newblock {Stochastic smoothed profile method for modeling random roughness in
  flow problems}.
\newblock {\em Computer Methods in Applied Mechanics and Engineering},
  263:99--112, 2013.

\end{thebibliography}
\end{document}